\documentclass[10pt]{amsart}

\usepackage{latexsym}
\usepackage{amsmath}
\usepackage{amsfonts}
\usepackage{amssymb}
\usepackage{esint}
\renewcommand{\a }{\alpha }
\renewcommand{\b }{\beta }
\renewcommand{\d}{\delta }
\newcommand{\D }{\Delta }

\newcommand{\e }{\varepsilon }
\newcommand{\g }{\gamma}

\renewcommand{\l }{\lambda }
\renewcommand{\L }{\Lambda }

\newcommand{\n }{\nabla }
\newcommand{\var }{\varphi }

\newcommand{\s }{\sigma }
\newcommand{\Sg }{\Sigma}
\renewcommand{\t }{\tau }

\renewcommand{\O }{\Omega }

\newcommand{\ov}{\overline}

\newcommand{\be}{\begin{equation}}
\newcommand{\ee}{\end{equation}}
\newenvironment{pf}{\noindent{\sc Proof}.\enspace}{\rule{2mm}{2mm}\medskip}
\newenvironment{pfn}{\noindent{\sc Proof}}{\rule{2mm}{2mm}\medskip}

\newcommand{\R}{\mathbb{R}}

\newcommand{\N}{\mathbb{N}}

\newcommand{\dis}{\displaystyle}
\newcommand{\Sf}{\mathbb{S}}
\newcommand{\dotcup}{\ \ensuremath{\mathaccent\cdot\cup}\,}
\newcommand{\T}{\mathbb{T}}

\begin{document}

\newtheorem{lem}{Lemma}[section]
\newtheorem{pro}[lem]{Proposition}
\newtheorem{thm}[lem]{Theorem}
\newtheorem{rem}[lem]{Remark}
\newtheorem{cor}[lem]{Corollary}
\newtheorem{df}[lem]{Definition}

\title[The Toda System on Compact Surfaces] {A variational Analysis of the Toda System on Compact Surfaces}

\author{Andrea Malchiodi and David Ruiz}

\address{SISSA, via Bonomea 265, 34136 Trieste (Italy) and Departamento de
An\'alisis Matem\'atico, University of Granada, 18071 Granada
(Spain).}

\thanks{A. M. has been partially supported by GENIL (SPR) for a stay in Granada in 2011,
and is supported by the FIRB project {\em Analysis and Beyond}
from MIUR. D.R has been supported by the Spanish Ministry of
Science and Innovation under Grant MTM2008-00988 and by J.
Andalucia (FQM 116).}

\email{malchiod@sissa.it, daruiz@ugr.es}



\begin{abstract}
In this paper we consider the following {\em Toda system} of
equations on a compact surface:
$$ \left\{
    \begin{array}{ll}
      - \D u_1 = 2 \rho_1 \left( \frac{h_1 e^{u_1}}{\int_\Sg
    h_1 e^{u_1} dV_g} - 1 \right) - \rho_2 \left( \frac{h_2 e^{u_2}}{\int_\Sg
    h_2 e^{u_2} dV_g} - 1 \right), \\
     - \D u_2 = 2 \rho_2 \left( \frac{h_2 e^{u_2}}{\int_\Sg
    h_2 e^{u_2} dV_g} - 1 \right) - \rho_1 \left( \frac{h_1 e^{u_1}}{\int_\Sg
    h_1 e^{u_1} dV_g} - 1 \right). &
    \end{array}
  \right.$$
We will give existence results by using variational methods in a
non coercive case. A key tool in our analysis is a new
Moser-Trudinger type inequality under suitable conditions on the
center of mass and the scale of concentration of the two
components $u_1, u_2$.
\end{abstract}

\maketitle

\section{Introduction}

Let $\Sg$ be a compact orientable surface without boundary, and
$g$ a Riemannian metric on $\Sg$. Consider the following system of
equations:
\begin{equation}\label{eq:toda}
    - \frac 1 2 \D u_i(x) = \sum_{j=1}^{N} a_{ij} e^{u_j(x)}, \qquad x \in \Sg, \ i = 1, \dots,
    N,
\end{equation}
where $\D=\D_g$ stands for the Laplace-Beltrami operator and $A =
(a_{ij})_{ij}$ is the {\em Cartan matrix} of $SU(N+1)$,
$$
  A = \left(
        \begin{array}{cccccc}
          2 & -1 & 0 & \dots & \dots & 0 \\
          -1 & 2 & -1 & 0 & \dots & 0 \\
          0 & -1 & 2 & -1 & \dots & 0 \\
          \dots & \dots & \dots & \dots & \dots & \dots \\
          0 & \dots & \dots & -1 & 2 & -1 \\
          0 & \dots & \dots & 0 & -1 & 2 \\
        \end{array}
      \right).
$$
Equation \eqref{eq:toda} is known as the {\em Toda system}, and
has been extensively studied in the literature. This problem has a
close relationship with geometry, since it can be seen as the
Frenet frame of holomorphic curves in $\mathbb{CP}^N$ (see
\cite{guest}). Moreover, it arises in the study of the non-abelian
Chern-Simons theory in the self-dual case, when a scalar Higgs field
is coupled to a gauge potential, see \cite{dunne, tar, yys}.

Let us assume, for the sake of simplicity, that $\Sg$ has total
area equal to $1$, i.e. $\int_{\Sg} 1 \, dV_g=1$. In this paper we
study the following version of the Toda system for $N=2$:


\begin{equation}\label{eq:gtodaaux}
  \left\{
    \begin{array}{l}
      - \D u_1 = 2 \rho_1 \left( h_1 e^{u_1} - 1 \right) - \rho_2 \left( h_2 e^{u_2} - 1 \right),   \\
     - \D u_2 = 2 \rho_2 \left(h_2 e^{u_2} - 1 \right) - \rho_1 \left( h_1 e^{u_1} - 1 \right),
    \end{array}
  \right. \end{equation}
where $h_i$ are smooth and strictly positive functions defined on
$\Sg$. By integrating on $\Sg$ both equations, we obtain that any
solution $(u_1,u_2)$ of \eqref{eq:gtodaaux} satisfies:
$$ \int_{\Sg} h_i e^{u_i} \, dV_g =1, \qquad i=1,\ 2. $$
Hence, problem \eqref{eq:gtodaaux} is equivalent to:

\begin{equation}\label{eq:gtoda}
  \left\{
    \begin{array}{ll}
      - \D u_1 = 2 \rho_1 \left( \frac{h_1 e^{u_1}}{\int_\Sg
    h_1 e^{u_1} dV_g} - 1 \right) - \rho_2 \left( \frac{h_2 e^{u_2}}{\int_\Sg
    h_2 e^{u_2} dV_g} - 1 \right), \\
     - \D u_2 = 2 \rho_2 \left( \frac{h_2 e^{u_2}}{\int_\Sg
    h_2 e^{u_2} dV_g} - 1 \right) - \rho_1 \left( \frac{h_1 e^{u_1}}{\int_\Sg
    h_1 e^{u_1} dV_g} - 1 \right). &
    \end{array}
  \right.
\end{equation}

Problem \eqref{eq:gtoda} is variational, and solutions can be
found as critical points of a functional $J_\rho : H^1(\Sg) \times
H^1(\Sg) \to \R$ ($\rho=(\rho_1,\rho_2)$) given by
\begin{equation} \label{funzionale}
J_\rho(u_1, u_2) = \int_\Sg Q(u_1,u_2)\, dV_g  + \sum_{i=1}^2
\rho_i \left ( \int_\Sg u_i dV_g - \log \int_\Sg h_i e^{u_i} dV_g
\right ),
\end{equation}
 where $Q(u_1,u_2)$ is defined as:
\begin{equation}\label{eq:QQ}
Q(u_1,u_2) = \frac{1}{3} \left ( |\n u_1|^2 + |\n u_2|^2 + \n u_1 \cdot
\n u_2\right ).
\end{equation}

Here and throughout the paper $\n u= \n_g u$ stands for the
gradient of $u$ with respect to the metric $g$, whereas $\cdot$
denotes the Riemannian scalar product.

Observe that both \eqref{eq:gtoda} and \eqref{funzionale} are
invariant under addition of constants to $u_1$, $u_2$. The
structure of the functional $J_{\rho}$  strongly depends on the
parameters $\rho_1$, $\rho_2$. To start with, the following
analogue of the Moser-Trudinger inequality has been given in
\cite{jw}:

\begin{equation} \label{mtjw}  4\pi \sum_{i=1}^2 \left (\log \int_\Sg h_i e^{u_i} dV_g -
\int_\Sg u_i dV_g \right ) \leq \int_\Sg Q(u_1,u_2)\, dV_g +C,
\end{equation}
for some $C=C(\Sg)$. As a consequence, $J_{\rho}$ is bounded from
below for $\rho_i \leq 4 \pi$ (see also \cite{sw, sw2, wang} for
related inequalities). In particular, if $\rho_i < 4 \pi$
($i=1,2$), $J_{\rho}$ is coercive and a solution for
\eqref{eq:gtoda} can be easily found as a minimizer.

If $\rho_i>4\pi$ for some $i=1,\ 2$, then $J_{\rho}$ is unbounded
from below and a minimization technique is no more possible. Let
us point out that the Leray-Schauder degree associated to
\eqref{eq:gtoda} is not known yet. For the scalar case, the
Leray-Schauder has been computed in \cite{clin}. The unique result
on the topological degree for Liouville systems is
\cite{lin-zhang}, but our case is not covered there. In this paper
we use variational methods to obtain existence of critical points
(generally of saddle type) for $J_{\rho}$.

\medskip Before stating our results, let us comment briefly on some aspects of the
problem under consideration. When some of
the parameters $\rho_i$ equals $4 \pi$, the situation becomes more
subtle. For instance, if we fix $\rho_1 < 4 \pi$ and let $\rho_2
\nearrow 4 \pi$, then $u_2$ could exhibit a blow-up behavior (see
the proof of Theorem 1.1 in \cite{jlw}). In this case, $u_2$ would
become close to a function $U_{\lambda, x}$ defined as:

$$ U_{\lambda, x}(y)= \log \left( \frac{4 \l}{\left(1 + \l \, d(x,y)^2\right)^2}
  \right),  $$
where $y \in \Sg$, $d(x,y)$ stands for the geodesic distance and
$\lambda$ is a large parameter. Those functions $U_{\l, x}$ are
the unique entire solutions of the Liouville equation (see
\cite{cl2}):

$$ - \D U= 2 e^{U}, \qquad \int_{\R^2} e^{U} \, dx <
+\infty.$$

In \cite{jlw} and \cite{ll} some conditions for existence are
given when some of the $\rho_i$'s equals $4 \pi$. The proofs involve a delicate analysis of the
limit behavior of the solutions when $\rho_i$ converge to $4\pi$
from below, in order to avoid bubbling of solutions. For that,
some conditions on the functions $h_i$ are needed.

The scalar counterpart of \eqref{eq:gtoda} is a Liouville-type problem in the form:
\begin{equation} \label{scalar} - \Delta u = 2 \rho \left( \frac{h(x) e^{u}}{\int_\Sigma h(x)
e^{u} d V_g} - 1 \right),\end{equation} with $\rho \in \R$. This
equation has been very much studied in the literature; there are
by now many results regarding existence, compactness of solutions,
bubbling behavior, etc. We refer the interested reader to the
reviews \cite{mreview, tar3}.

Solutions of \eqref{scalar} correspond to critical points of the
functional $I_{\rho}:H^1(\Sg) \to \R$,
\begin{equation}\label{scalar2}
    I_\rho(u) =  \frac 1 2 \int_{\Sigma} |\n_g u|^2 dV_g + 2 \rho \left ( \int_\Sigma
    u dV_g -  \log \int_\Sigma h(x) e^{u} dV_g \right ),
    \qquad u \in H^{1}(\Sigma).
\end{equation}

The classical Moser-Trudinger inequality implies that $I_{\rho}$
is bounded from below for $\rho \leq 4 \pi$. For larger values of
$\rho$, variational methods were applied to \eqref{scalar} for the
first time in \cite{djlw}, \cite{st}. In \cite{dm} the
$Q$-curvature prescription problem is addressed in a 4-dimensional
compact manifold: however, the arguments of the proof can be
easily translated to the Liouville problem \eqref{scalar}, see
\cite{dja}.

Let us  briefly describe the proof of \cite{dm} in the case $\rho
\in (4 \pi , 8 \pi)$, for simplicity. In \cite{dm} it is shown
that, whenever $I_{\rho}(u_n) \to -\infty$, then (up to a
subsequence)
$$ \frac{e^{u_n}}{\int_{\Sg} e^{u_n}\, dV_g} \rightharpoonup
\delta_{x}, \ x \in \Sg,
$$ in the sense of measures. Moreover, for $L>0$ sufficiently large, one can define a
homotopy equivalence (see also \cite{mal}):

$$I_{\rho}^{-L}= \{u\in H^1(\Sg):\ I_{\rho}(u)< -L \} \simeq \{ \delta_x:\ x \in \Sg\} \simeq \Sg.$$

Therefore the sublevel $I_{\rho}^{-L}$ is not contractible, and
this allows us to use a min-max argument to find a solution. We
point out that \cite{dm} also deals with the case of higher values
of $\rho$, whenever $\rho \notin 4 \pi \N$.


\medskip Coming back to  system \eqref{eq:gtoda}, there are very
few results when $\rho_i>4\pi$ for some $i=1$, $2$. One of them is
given in \cite{cheikh} and concerns the case $\rho_1<4\pi$ and
$\rho_2 \in (4\pi m, 4 \pi (m+1))$, $m \in \N$. There, the
situation is similar to \cite{dm}; in a certain sense, one can
describe the set $J_{\rho}^{-L}$ from the behavior of the second
component $u_2$ as in \cite{dm}.

In Theorem 1.4 of \cite{jlw}, an existence result is stated for
$\rho_i \in (0,4 \pi) \cup (4 \pi, 8 \pi)$ for a compact surface
$\Sg$ with positive genus: however, the min-max argument used in
the proof seems not to be correct. The main problem is that a
one-dimensional linking argument is used to obtain conditions on
both the components of the system. In any case, the core of
\cite{jlw} is the blow-up analysis for the Toda system (see Remark
\ref{puffff} for more details). In particular, it is shown that if
the $\rho_i$'s are bounded away from $4 \pi \N$, the set of
solutions of \eqref{eq:gtoda} is compact (up to addition of
constants). This is an essential tool for our analysis.

In this paper we deal with the case $\rho_i \in (4 \pi, 8 \pi)$,
$i=1,2$. Our main result is the following:

\begin{thm}\label{t:main} Assume that $\rho_i \in (4 \pi, 8 \pi)$ and that $h_1, h_2$
are two positive $C^1$ functions on $\Sg$. Then there exists a
solution $(u_1, u_2)$ of \eqref{eq:gtoda}.
\end{thm}

Let us point out that we find existence of solutions also if $\Sg$
is a sphere. Moreover, our existence result is based on a detailed
study of the topological properties of the low sublevels of
$J_{\rho}$. This study is interesting in itself; in the scalar
case an analogous one has been used to deduce multiplicity
results (see \cite{demarchis}) and degree computation formulas
(see \cite{mal}).

We shall see that the low sublevels of $J_{\rho}$ contain couples
in which at least one component is very concentrated around some
point of $\Sg$. Moreover, both components can concentrate at two
points that could eventually coincide. However, we shall see that,
in a certain sense,
\begin{equation} \label{frase} \mbox{\em if } u_1,\ u_2 \ \mbox{\em concentrate around the same point at
the same rate, then } J_{\rho} \, \mbox{\em is bounded from below.}
\end{equation}

To make this statement rigorous, we need several tools.

The first is a definition of a rate of concentration of a positive
function $f \in \Sg$, normalized in $L^1$, which is a refinement
of the one given in \cite{mr}; this will be measured by a positive
parameter called $\s=\s(f)$. In a sense, the smaller is $\s$, the
higher is the rate of concentration of $f$. Compared to the
classical concentration compactness arguments, our function $\s$
has the property of being continuous with respect to the $L^1$
topology (see Remark \ref{sigmacont}). Second, we also need to
define a continuous center of mass when $\s\leq \d$ for some fixed
$\d>0$: we will denote it by $\beta=\beta(f) \in \Sg$. When $\s
\geq \d$, the function is not concentrated and the center of mass
cannot be defined. Hence, we have a map:

$$\psi: H^1(\Sg) \to \ov{\Sg}_{\d}, \ \psi(u_i)= (\beta(f_i), \s(f_i)),\ \mbox{ where } f_i=\frac{e^{u_i}}{\int_{\Sg}
e^{u_i}\, dV_g}.$$ Here $\ov{\Sg}_{\d}$ is the topological cone
with base $\Sg$, so that we make the identification to a point
when $\s \geq \delta $ for some $\delta>0$ fixed.

Third, we need an improvement of the Moser-Trudinger inequality in
the following form: if $ \psi(f_1) = \psi(f_2)$, then
$J_{\rho}(u_1,u_2)$ is bounded from below. In this sense,
\eqref{frase} is made precise. The proof uses local versions of
the Moser-Trudinger inequality and applications of it to small
balls (via a convenient dilation) and to annuli with small
internal radius (via a Kelvin transform).

Roughly speaking, on  low sublevels one of the following
alternatives hold:

\begin{enumerate}

\item one component concentrates at a point whereas the other does
not concentrate ($\s_i< \delta \leq \s_j$), or

\item the two components concentrate at different points ($\s_i <
\delta,\ \beta_1 \neq \beta_2$), or

\item the two components concentrate at the same point with
different rates of concentration ($\s_i< \s_j<\delta$,
$\beta_1=\beta_2$).

\end{enumerate}

With this at hand, for $L > 0$ large we are able to define a  continuous
map:
$$ J_{\rho}^{-L} \quad \stackrel{\psi \oplus \psi}{\longrightarrow}
\quad X:=(\ov{\Sg}_{\d} \times \ov{\Sg}_{\d}) \setminus \ov{D},
$$
where $\ov{D}$ is the diagonal of $\ov{\Sg}_{\d} \times
\ov{\Sg}_{\d}$. We can also proceed in the opposite direction: in
Section \ref{s:4} we construct a family of test functions modeled
on $X$ on which $J_\rho$ attains arbitrarily low values, see Lemma
\ref{l:dsmallIlow} for the precise result. Calling $\phi : X \to
J_{\rho}^{-L}$ the corresponding map, we will prove that the
composition
\begin{equation} \label{compo}
X \quad \stackrel{\phi}{\longrightarrow} \quad J_{\rho}^{-L} \quad
\stackrel{\psi \oplus \psi}{\longrightarrow} \quad X
\end{equation}
is homotopically equivalent to the identity map. In this situation
it is said that $J_{\rho}^{-L}$ {\em dominates} $X$ (see
\cite{hat}, page 528). In a certain sense, those maps are natural
since they describe properly the topological properties of
$J_{\rho}^{-L}$.

We will see that for any compact orientable surface $\Sg$, $X$ is
non-contractible; this is proved by estimating its cohomology
groups. As a consequence, $\phi(X)$ is not contractible in
$J_{\rho}^{-L}$. This allows us to use a min-max argument to find
a critical point of $J_{\rho}$. Here, the compactness of solutions
proved in \cite{jlw} is an essential tool, since the Palais-Smale
property for $J_{\rho}$ is an open problem (as it is for the
scalar case).

The rest of the paper is organized as follows. In Section
\ref{s:pr} we present the notations that will be used in the
paper, as well as some preliminary results. The definition of the
map $\psi$, its properties, and the improvement of the
Moser-Trudinger inequality will be exposed in Section \ref{s:3}.
In Section \ref{s:4} we define the map $\phi$ and prove that the
composition \eqref{compo} is homotopic to the identity. Here we
also develop the min-max scheme that gives a critical point of
$J_{\rho}$. The fact that $X$ is not contractible is proved in a
final Appendix.

\section{Notations and preliminaries}\label{s:pr}

In this section we collect some useful notation and preliminary
facts. Throughout the paper, $\Sg$ is a compact orientable surface
without boundary; for simplicity, we assume $|\Sg|= \int_{\Sg} 1
dV_g =1$. Given $\delta>0$, we define the topological cone:
\begin{equation} \label{cono} \ov{\Sg}_{\delta} = \left(\Sg \times (0, +\infty) \right)
|_{\left( \Sg \times [\delta, + \infty) \right)}.\end{equation}

For $x, y \in \Sg$ we denote by $d(x,y)$ the metric distance
between $x$ and $y$ on $\Sg$. In the same way, for any $p \in \Sg
$, $\Omega, \Omega' \subseteq \Sg$, we denote:

$$ d(p, \Omega) = \inf \left\{ d(p,x) \; :  x \in \Omega
\right\}, \qquad
  d(\Omega,\Omega') = \inf \left\{ d(x,y) \; : \; x \in \Omega,\ y \in
  \Omega'
\right\}.
$$

Moreover, the symbol $B_p(r)$ stands for the open metric ball of
radius $r$ and center $p$, and $A_p(r,R)$ the open annulus of
radii $r$ and $R$, $r<R$. The complement of a set $\Omega$ in
$\Sg$ will be denoted by $\Omega^c$.

Given a function $u \in L^1(\Sg)$ and $\Omega \subset \Sg$, we
consider the average of $u$ on $\Omega$:
$$ \fint_{\Omega} u \, dV_g = \frac{1}{|\Omega|} \int_{\Omega} u \,
dV_g.$$ We denote by $\ov{u}$ the average of $u$ in $\Sg$: since
we are assuming $|\Sg| = 1$, we have
$$
 \ov{u}= \int_\Sg u \, dV_g = \fint_\Sg u\, dV_g.
$$

Throughout the paper we will denote by $C$ large constants which
are allowed to vary among different formulas or even within lines.
When we want to stress the dependence of the constants on some
parameter (or parameters), we add subscripts to $C$, as $C_\d$,
etc.. Also constants with subscripts are allowed to vary.
Moreover, sometimes we will write $o_{\alpha}(1)$ to denote
quantities that tend to $0$ as $\alpha \to 0$ or $\alpha \to
+\infty$, depending on the case. We will similarly use the symbol
$O_\a(1)$ for bounded quantities.

\

\noindent We begin by recalling the following compactness result from \cite{jlw}.

\begin{thm}\label{th:jlw} (\cite{jlw})
Let $m_1, m_2$ be two non-negative integers, and suppose $\L_1,
\L_2$ are two compact sets of the intervals $(4 \pi m_1, 4 \pi (m_1
+ 1))$ and $(4 \pi m_2, 4 \pi (m_2 + 1))$ respectively. Then if
$\rho_1 \in \L_1$ and $\rho_2 \in \L_2$ and if we impose
$\int_{\Sg} u_i dV_g = 0$, $i = 1, 2$, the solutions of
\eqref{eq:gtoda} stay uniformly bounded in $L^\infty(\Sg)$
(actually in every $C^l(\Sg)$ with $l \in \N$).
\end{thm}

\noindent Next, we also recall some Moser-Trudinger type
inequalities. As commented in the introduction, problem
\eqref{eq:gtoda} is the Euler-Lagrange equation of the energy
functional $J_{\rho}$ given in \eqref{funzionale}. This functional
is bounded below only for certain values of $\rho_1, \rho_2$, as
has been proved by Jost and Wang (see also \eqref{mtjw}):

\begin{thm}\label{th:jw} (\cite{jw})
The functional $J_\rho$ is bounded from below if and only if
$\rho_i \leq 4 \pi$, $i=1,\ 2$.
\end{thm}


\noindent The next proposition can be thought of as a local
version of Theorem \ref{th:jw}, and will be of use in Section
\ref{s:3}. Let us recall the definition of the quadratic form $Q$ in \eqref{eq:QQ}.

\begin{pro}\label{p:MTbd}  Fix $\d > 0$, and let $\Omega_1 \subset \Omega_2 \subset \Sg$ be such that
$d(\Omega_1, \partial \Omega_2) \geq \d$. Then, for any $\e > 0$
there exists a constant $C = C(\e, \delta)$ such that for all $u
\in H^1(\Sg)$
\begin{equation}\label{eq:ineqSg}
4 \pi \left ( \log \int_{\Omega_1}  e^{u_1} dV_g + \log
\int_{\Omega_1} e^{u_2} dV_g -\fint_{\Omega_2} u_1 dV_g -
\fint_{\Omega_2} u_2 dV_g \right )\leq (1+\e) \int_{\Omega_2}
Q(u_1, u_2) dV_g + C.
\end{equation}
\end{pro}

\begin{pf}
We can assume without loss of generality that $\fint_{\Omega_2}
u_i dV_g = 0$ for $i = 1, 2$. Let us write
$$
u_i = v_i + w_i, \quad \int_{\Omega_2} v_i \, dV_g =
\int_{\Omega_2} w_i \, dV_g =0,
$$
where $v_i \in L^\infty(\Omega_2)$ and $w_i \in H^1(\Omega_2)$
will be fixed later. We have
\begin{equation}\label{eq:ddmm2}
    \log \int_{\Omega_1} e^{u_1} dV_g + \log \int_{\Omega_1} e^{u_2} dV_g \leq
    \| v_1\|_{L^{\infty}(\Omega_1)} + \| v_2\|_{L^{\infty}(\Omega_1)}
    + \log \int_{\Omega_1} e^{w_1} dV_g + \log \int_{\Omega_1} e^{w_2} dV_g.
\end{equation}
We next consider a smooth cutoff function $\chi$ with values into
$[0,1]$ satisfying
$$
  \left\{
    \begin{array}{ll}
      \chi(x) = 1 & \hbox{ for } x \in \Omega_1,\\
      \chi(x) = 0 & \hbox{ if } d(x, \Omega) > \delta/2,
    \end{array}
  \right.
$$
and then define
$$
  \tilde{w}_i(x) = \chi(x) w_i(x); \qquad \quad i = 1, 2.
$$
Clearly $\tilde{w}_i$ belongs to $H^1(\Sg)$ and is supported in a
compact set of the interior of $\Omega_2$. Hence we can apply
Theorem \ref{th:jw} to $\tilde{w}_i$ on $\Sg$, finding
$$
 \log \int_{\Omega_1} e^{w_1} dV_g + \log \int_{\Omega_1} e^{w_2} dV_g \leq \log
 \int_{{\Sg}} e^{\tilde{w}_1} dV_g + \log \int_{{\Sg}} e^{\tilde{w}_2} dV_g
 \leq $$$$\frac{1}{4\pi} \int_{{\Sg}} Q(\tilde{w}_1,\tilde{w}_2) dV_g + \fint_{\Sg} (\tilde{w}_1 + \tilde{w}_2)\, dV_g + C.
$$
Using the Leibnitz rule and H\"older's inequality we obtain
$$
  \int_{{\Sg}} Q(\tilde{w}_1,\tilde{w}_2) dV_g \leq (1+\e) \int_{\Omega_2}
  Q(w_1, w_2) dV_g + C_{\e} \int_{\Omega_2}
  (w_{1}^2 + w_{2}^2) dV_g.
$$

Moreover, we can estimate the mean value of $\tilde{w}_i$ in the
following way:

$$\fint_{\Sg} \tilde{w}_i \, dV_g \leq  C \left ( \int_{\Sg} |\n \tilde{w}_i|^2\, dV_g \right )^{1/2} \leq
C_{\e} + \e \int_{\Omega_2} |\n \tilde{w}_i|^2\, dV_g \leq
$$$$C_{\e} + C \e \left ( \int_{\Omega_2} |\n w_i|^2\, dV_g + C
\int_{\Omega_2} w_{i}^2\, dV_g \right ).$$

From \eqref{eq:ddmm2} and the last formulas we find
\begin{eqnarray}\label{eq:last} \nonumber
  \log \int_{\Omega_1}  e^{u_1} dV_g + \log \int_{\Omega_1}  e^{u_2} dV_g &
  \leq & \| v_1\|_{L^{\infty}(\Omega_1)} + \| v_2\|_{L^{\infty}(\Omega_1)}
  + \frac{1+\e}{4\pi} \int_{\Omega_2} Q(w_1, w_2) dV_g \\ & + &  C_{\e}
  \int_{\Omega_2} (w_1^2 + w_2^2) dV_g + C.
\end{eqnarray}
To control the latter terms we use truncations in Fourier modes.
Define $V_{\e}$ to be the direct sum of the eigenspaces of the
Laplacian on $\Omega_2$ (with Neumann boundary conditions) with
eigenvalues less or equal than $C_{\e}\e^{-1}$. Take now $v_i$ to
be the orthogonal projection of $u_i$ onto $V_{\e}$. In $V_\e$ the
$L^\infty$ norm is equivalent to the $L^2$ norm: by using
Poincar{\'e}'s inequality we get
$$ C_{\e} \int_{\Omega_2} (w_{1}^2 + w_{2}^2) dV_g  \leq \e \int_{\Omega_2} Q(u_1,u_2) dV_g, $$$$
  \|v_i\|_{L^\infty(\Omega_1)} \leq C_{\e} \|v_i\|_{L^2(\Omega_2)} \leq C_\e \left(
  \int_{\Omega_2} |\n u_i|^2 dV_g \right)^{\frac 12} \leq
\e \int_{\Omega_2} Q(u_1,u_2) dV_g +C_{\e}.
$$
Hence, from \eqref{eq:last} and the above inequalities we derive
\eqref{eq:ineqSg} by renaming $\e$ properly.

\end{pf}

\begin{rem}\label{r:regeigenv} While the Fourier decomposition used in the above proof
depends on  $\Omega_2$, the constants only depend on $\Sg$, $\d$
and $\e$. In fact, one can replace $\O_2$ by a domain $\check{\O}_2$,
$\O_2 \subseteq \check{\O}_2 \subseteq B_{\O_2}(\d/2)$ with boundary
curvature depending only on $\d$ and satisfying a uniform interior
sphere condition with spheres of radius $\d^3$. For example, one
can obtain such a domain $\check{\O}_2$ triangulating $\Sg$ by
simplexes with diameters of order $\d^2$, take suitable union of
triangles and smoothing the corners. For these domains, which are
finitely many, the eigenvalue estimates will only depend on $\d$.
\end{rem}

\

\noindent We next prove a criterion which gives us a first insight
on the properties of the low sublevels of $J_{\rho}$. This result
is in the spirit of an improved inequality in \cite{cl}, and we
use an extra covering argument to track the concentration
properties of both components of the system. We need first an
auxiliary lemma.

\begin{lem}\label{l:step1}
Let $\d_0 > 0$, $\g_0 > 0$, and let $\Omega_{i,j} \subseteq \Sg$,
$i,j = 1, 2$, satisfy $d(\Omega_{i,j},\Omega_{i,k}) \geq \d _0$
for $j \neq k$. Suppose that $u_1, u_2 \in H^1(\Sg)$ are two
functions verifying
\begin{equation}\label{eq:ddmm}
    \frac{\int_{\Omega_{i,j}} e^{u_i} dV_g}{\int_\Sg e^{u_i} dV_g}
    \geq \g_0, \qquad \qquad i,j = 1, 2.
\end{equation}
Then there exist positive constants $\tilde{\g}_0$, $\tilde{\d}_0$,
depending only on $\g_0$, $\d_0$, and two sets $\tilde{\O}_1, \tilde{\O}_2
\subseteq \Sg$, depending also on $u_1, u_2$ such that
\begin{equation}\label{eqLddmm2}
    d(\tilde{\O}_1, \tilde{\O}_2) \geq \tilde{\d}_0; \qquad
    \quad \frac{\int_{\tilde{\Omega}_{i}} e^{u_1} dV_g}{\int_\Sg e^{u_1} dV_g}
    \geq \tilde{\g}_0, \quad \frac{\int_{\tilde{\Omega}_{i}} e^{u_2} dV_g}{\int_\Sg
    e^{u_2} dV_g} \geq \tilde{\g}_0; \quad i = 1, 2.
\end{equation}
\end{lem}

\begin{pf}
First, we fix a number $r_0 < \frac{\d_0}{80}$. Then we cover
$\Sg$ with a finite union of metric balls $(B_{x_l}(r_0))_l$,
whose number can be bounded by an integer $N_{r_0}$ which depends
only on $r_0$ (and $\Sg$).

Next we cover $\ov{\Omega}_{i,j}$ by a finite number of these
balls, and we choose $y_{i,j} \in \cup_l (x_l)$ such that
$$
  \int_{B_{y_{i,j}}(r_0)} e^{u_i} dV_g = \max \left\{ \int_{B_{x_l}(r_0)}
  e^{u_i} dV_g \; : \; B_{x_l}(r_0) \cap \ov{\Omega}_{i,j} \neq \emptyset
  \right\}.
$$
Since the total number of balls is bounded by $N_{r_0}$ and since
by our assumption the (normalized) integral of $e^{u_i}$ over
$\Omega_{i,j}$ is greater or equal than $\g_0$, it follows that
\begin{equation}\label{eq:bdtg0}
  \frac{\int_{B_{y_{i,j}}(r_0)} e^{u_i} dV_g}{\int_{\Sg}
  e^{u_i} dV_g} \geq \frac{\g_0}{N_{r_0}}.
\end{equation}
By the properties of the sets $\Omega_{i,j}$, we have that:
$$
 B_{y_{i,j}}(2r_0) \cap B_{y_{i,k}}(r_0)=\emptyset \qquad
 \hbox{ for } j \neq k.
$$
Now, one of the following two possibilities occurs:

\begin{description}
  \item[(a)] $B_{y_{1,1}}(5 r_0) \cap \left( B_{y_{2,1}}(5 r_0)
  \cup B_{y_{2,2}}(5 r_0) \right) \neq \emptyset$ or $B_{y_{1,2}}(5 r_0) \cap \left( B_{y_{2,1}}(5 r_0)
  \cup B_{y_{2,2}}(5 r_0) \right) \neq \emptyset$;
  \item[(b)] $B_{y_{1,1}}(5 r_0) \cap \left( B_{y_{2,1}}(5 r_0)
  \cup B_{y_{2,2}}(5 r_0) \right) = \emptyset$ and $B_{y_{1,2}}(5 r_0) \cap \left( B_{y_{2,1}}(5 r_0)
  \cup B_{y_{2,2}}(5 r_0) \right) = \emptyset$.
\end{description}
In case {\bf (a)} we define the sets $\tilde{\Omega}_i$ as
$$
  \tilde{\Omega}_1 = B_{y_{1,1}}(30 r_0), \qquad \quad
  \tilde{\Omega}_2 = B_{y_{1,1}}(40 r_0)^c,
$$
while in case {\bf (b)} we define
$$
  \tilde{\Omega}_1 = B_{y_{1,1}}(r_0) \cup B_{y_{2,1}}(r_0);
  \qquad \quad \tilde{\Omega}_2 = B_{y_{1,2}}(r_0) \cup B_{y_{2,2}}(r_0)).
$$
We also set $\tilde{\g}_0 = \frac{\g_0}{N_{r_0}}$ and
$\tilde{\d}_0 = r_0$. We notice that $\tilde{\g}_0$ and
$\tilde{\d}_0$ depend only on $\g_0$ and $\d_0$, as claimed, and
that the sets $\tilde{\Omega}_i$ satisfy the required conditions.
This concludes the proof of the lemma.

\end{pf}

\

\noindent We next derive the improvement of the constants in Theorem
\ref{th:jw}, in the spirit of \cite{cl}.

\begin{pro}\label{p:imprc}

Let $u_1, u_2 \in H^1(\Sg)$ be a couple of functions satisfying
the assumptions of Lemma \ref{l:step1} for some positive constants
$\d_0, \g_0$. Then for any $\e > 0$ there exists $C=C(\e) > 0$,
depending on $\e, \d_0$, and $\g_0$ such that
$$
  8 \pi \left( \log \int_{\Sg} e^{u_1-\ov{u}_1} dV_g + \log
  \int_{\Sg} e^{u_2-\ov{u}_2} dV_g \right) \leq (1+\e) \int_\Sg
  Q(u_1, u_2) dV_g + C.
$$
\end{pro}

\begin{pf} Let $\tilde{\d}_0, \tilde{\g}_0$ and $\tilde{\Omega}_1, \tilde{\Omega}_2$ be as
in Lemma \ref{l:step1}, and assume without loss of generality that
$\ov{u}_1 = \ov{u}_2 = 0$. Let us define $U_i=\{x \in \Omega:\
d(x, \tilde{\Omega}_i) < \tilde{\d}_0/2\}$. By applying
Proposition \ref{p:MTbd}, we get:

\begin{equation} \label{hola} 4 \pi \left (\log \int_{\tilde{\Omega}_i}  e^{u_1} dV_g + \log
\int_{\tilde{\Omega}_i} e^{u_2} dV_g - \fint_{U_i} (u_1+u_2) dV_g
\right)\leq (1+\e) \int_{U_i} Q(u_1, u_2) dV_g + C.\end{equation}

Observe that:
$$ \log \int_{\tilde{\Omega}_i}  e^{u_j} dV_g\geq \log \left ( \int_{\Sg}  e^{u_j}
dV_g\right ) + \log \tilde{\g}_0.$$

Since $U_1 \cap U_2 = \emptyset$, we can sum \eqref{hola} for
$i=1,2$, to obtain

$$8 \pi  \left ( \log \int_{\Sg}  e^{u_1} dV_g + \log
\int_{\Sg} e^{u_2} dV_g -\sum_{i=1}^2 \fint_{U_i} (u_1+u_2) dV_g
\right )\leq (1+\e) \int_{\Sg} Q(u_1, u_2) dV_g + C.$$

It suffices now to estimate the term $\fint_{U_i} (u_1+u_2) dV_g$.
By using Poincar{\'e}'s inequality and the estimate $|U_i| \geq
\tilde{\d}_0^2$, we have:
$$ \fint_{U_i} u_j \, dV_g \leq  \tilde{\d}_0^{-2} \int_{U_i} u_j \,
dV_g \leq C \left ( \int_{\Sg} |\n u_j|^2\, dV_g \right )^{1/2}
\leq C + \e \int_{\Sg} |\n u_j|^2\, dV_g.$$

To finish the proof it suffices to properly rename $\e$.

\end{pf}

Proposition \ref{p:imprc} implies that on low sublevels, at
least one of the components must be very concentrated around a
certain point. A more precise description of the topological
properties of $J_{\rho}^{-L}$ will be given later on.

\section{Volume concentration and improved inequality} \label{s:3}

In this section we give the main tools for the description of the
sublevels of the energy functional $J_{\rho}$. Those will be
contained in Propositions \ref{covering} and \ref{mt}, whose proof
will be given in the subsequent subsections.

First, we give continuous definitions of center of mass and scale
of concentration of positive functions normalized in $L^1$, which
are adequate for our purposes. Those are a refinement of
\cite{mr}.

Consider the set
$$
  A = \left\{ f \in L^1(\Sg) \; : \; f >
  0 \ \hbox{ a. e. and } \int_\Sg f dV_g = 1 \right\},
$$
endowed with the topology inherited from $L^1(\Sg)$.

Moreover, let us recall the definition \eqref{cono} for the cone
$\ov{\Sg}_{\delta}$.

\begin{pro} \label{covering} Let us fix a constant $R>1$. Then
there exists $\delta= \delta(R) >0$ and a continuous map:
$$ \psi : A \to \ov{\Sg}_{\delta}, \qquad \quad \psi(f)= (\beta, \sigma),$$
satisfying the following property: for any $f \in A$ there exists
$p \in \Sg$ such that

\begin{enumerate}

\item[{\emph a)}] $ d(p, \beta) \leq  C' \sigma$ for $C' = \max\{3
R+1, \delta^{-1}diam(\Sg) \}.$

\item[{\emph b)}] There holds: $$ \int_{B_p(\sigma)} f \, dV_g >
\tau, \qquad \quad \int_{B_p(R \sigma)^c} f \, dV_g > \tau, $$

where $\tau>0$ depends only on $R$ and $\Sg$.

\end{enumerate}

\end{pro}

Roughly speaking, the above map $\psi(f)= (\beta, \sigma)$ gives
us a center of mass of $f$ and its scale of concentration around
that point. Indeed, the smaller is $\sigma$, the bigger is the
rate of concentration. Moreover, if $\sigma$ exceeds a certain
positive constant, $\beta$ could not be defined; so, it is natural
to make the identification in $\ov{\Sg}_\d$.

Next, we state an improved Moser-Trudinger inequality for couples
$(u_1,u_2)$ such that $e^{u_i}$ are centered at the same point
with the same rate of concentration. Being more specific, we have
the following:

\begin{pro} \label{mt} Given any $\e>0$, there exist $R=R(\e)>1$ and $\psi$ as given in Proposition \ref{covering},
such that for any $(u_1, u_2) \in H^1(\Sg)\times H^1(\Sg)$ with:

$$\psi \left( \frac{e^{u_1}}{\int_{\Sigma} e^{u_1} dV_g} \right
)= \psi \left( \frac{e^{u_2}}{\int_{\Sigma} e^{u_2} dV_g} \right
),
$$

the following inequality holds:

$$ (1+\e) \int_\Sg Q(u_1,u_2)\, dV_g
\geq  8 \pi \left (\log \int_\Sg e^{u_1-\ov{u}_1} + dV_g +\log
\int_\Sg e^{u_2-\ov{u}_2} dV_g \right )+ C,
$$
for some $C=C(\e)$.
\end{pro}

The rest of the section is devoted to the proof of those
propositions.

\subsection{Proof of Proposition \ref{covering}}

Take $R_0=3R$, and define $ \sigma: A \times \Sg \to (0,+\infty)$
such that:

\begin{equation} \label{sigmax} \int_{B_x(\sigma(x,f))} f \, dV_g = \int_{B_x(R_0
\sigma(x,f))^c} f \, dV_g.  \end{equation}

It is easy to check that $\sigma(x,f)$ is uniquely determined and
continuous. Moreover, $\sigma$ satisfies:

\begin{equation} \label{dett} d(x,y) \leq R_0 \max \{ \sigma(x,f), \sigma(y,f)\}
+\min \{ \sigma(x,f), \sigma(y,f)\}.
\end{equation}
Otherwise, $ B_x(R_0 \sigma(x,f)) \cap B_y(\sigma(y,f)+\e) =
\emptyset $ for some $\e>0$. Moreover, since $B_y(R_0 \sigma(y,f))$
does not fulfil the whole space $\Sg$, $A_y(\sigma(y,f),
\sigma(y,f)+\e)$ is a nonempty open set. Then:
$$ \int_{B_x(\sigma(x,f))} f \, dV_g =
\int_{B_x(R_0 \sigma(x,f))^c} f \, dV_g \geq
\int_{B_y(\sigma(y,f)+\e)} f \, dV_g
> \int_{B_y(\sigma(y,f))} f \, dV_g. $$ By interchanging the roles of $x$ and $y$, we would also
obtain the reverse inequality. This contradiction proves
\eqref{dett}.

We now define:

$$
  T: A \times \Sg \to \R, \qquad T(x,f) = \int_{B_x(\s(x,f))} f dV_g.
$$

\begin{lem}\label{sigma} If $x_0 \in \Sg$ is such that $T(x_0,f) =
\max_{y \in\Sg} T(y,f)$, then $\s(x_0,f) < 3\s(x, f)$ for any
other $x \in \Sg$.
\end{lem}

\begin{pf} Choose any $x \in \Sg$ and $\e>0$. First, observe that $B_x(R_0 \s(x,f)+\e)
$ must intersect $B_{x_0}(\s(x_0,f))$. Otherwise, as above, we
know that $A_x(R_0 \s(x,f), R_0 \s(x,f)+\e)$ is an open nonempty
set. Then
$$ T(x_0,f)= \int_{B_{x_0}(\s(x_0,f))} f \, dV_g < \int_{B_{x}(R_0\s(x,f))^c} f \, dV_g =T(x,f),  $$
a contradiction. Arguing in the same way, we can also conclude
that $B_x(R_0 \s(x,f)+\e) $ cannot be contained in
$B_{x_0}(R_0\s(x_0,f))$.

By the triangular inequality, we obtain that:

$$ 2 (R_0 \sigma(x,f)+\e) > (R_0-1) \sigma(x_0,f).$$

Since $\e>0$ is arbitrary, there follows:

$$  \sigma(x,f) \geq \frac{R_0-1}{2 R_0} \sigma(x_0,f).$$

Recalling that $R_0>3$, we are done.

\end{pf}

As a consequence of the previous lemma, we obtain the following:

\begin{lem} \label{tau} There exists a fixed $\t > 0$ such that
$$
  \max_{x \in \Sg} T(x,f) > \t > 0 \qquad \quad \hbox{ for all } f \in A.
$$
\end{lem}

\begin{pf}
Let us fix $x_0 \in \Sigma$ such that $T(x_0,f)=\max_{x\in \Sigma}
T(x,f)$. For any $x \in A_{x_0}(\sigma(x_0,f ), R\sigma(x_0,f))$,
by Lemma \ref{sigma}, we have that:

$$ \int_{B_x(\sigma(x_0,f)/3)} f \, dV_g \leq \int_{B_x(\sigma(x,f))} f \, dV_g  \leq T(x_0,f). $$

Let us take a finite covering:

$$A_{x_0}(\sigma(x_0,f ), R\sigma(x_0,f)) \subset \cup_{i=1}^k
B_{x_i}(\sigma(x_0,f)/3).$$

Observe that $k$ is independent of $f$ or $\sigma(x_0,f)$, and
depends only on $\Sigma$ and $R$. Therefore:

$$ 1 = \int_{\Sg} f \, dV_g \leq \int_{B_{x_0}(\sigma(x_0,f))} f
\, dV_g + \int_{B_{x_0}(R\sigma(x_0,f))^c} f \, dV_g+ \sum_{i=1}^k
\int_{B_{x_i}(\sigma(x_0,f)/3)} f \, dV_g \leq (k+2) T(x_0,f).$$

\end{pf}

Let us define:

$$ \sigma : A \to \R, \qquad \quad \sigma(f)= 3 \min\{ \sigma(x,f): \ x
\in \Sg \},$$ which is obviously a continuous function.

\begin{rem} \label{sigmacont} In \cite{mr} (see Section 3 there) a sort
of concentration parameter is defined, but it does not depend continuously
on $f$. Moreover, the definition of barycenter given below has been modified
compared to \cite{mr}. Finally, the application
$\psi$ is mapped to a cone; this interpretation, which is crucial
in our framework, was missing in \cite{mr}.

\end{rem}

Given $\tau$ as in Lemma \ref{tau}, consider the set:
\begin{equation} \label{defS}
  S(f) = \left\{ x \in \Sg \; : \; T(x,f) > \t,\ \s(x,f) <  \s(f) \right\}.
\end{equation}

If $x_0 \in \Sg$ is such that $T(x_0,f)= \max_{x\in \Sg} T(x,f)$,
then Lemmas \ref{sigma} and \ref{tau} imply that $x_0 \in S(f)$.
Therefore, $S(f)$ is a nonempty open set for any $f \in A$.
Moreover, from \eqref{dett}, we have that:

\begin{equation}\label{eq:diamS} diam(S(f)) \leq (R_0+1)\s(f).\end{equation}

By the Nash embedding theorem, we can assume that $\Sg \subset
\R^N$ isometrically, $N \in \N$. Take an open tubular neighborhood $\Sg
\subset U \subset \R^N$ of $\Sg$, and $\delta>0$ small enough so
that:

\begin{equation} \label{co}
 co \left [ B_x((R_0+1)\delta)\cap \Sg \right ] \subset
 U \ \forall \, x \in \Sg, \end{equation}
where $co$ denotes the convex hull in $\R^N$.

We define now
$$
  \eta(f) = \frac{\displaystyle \int_\Sg (T(x,f) - \t)^+ \left( \s(f) - \s(x,f)
  \right)^+ x \ dV_g}{\displaystyle \int_\Sg (T(x,f) - \t)^+ \left( \s(f) - \s(x,f)
  \right)^+ dV_g}\in \R^N.
$$

The map $\eta$ yields a sort of center of mass in $\R^N$. Observe
that the integrands become nonzero only on the set $S(f)$.
However, whenever $\sigma(f) \leq \delta$, \eqref{eq:diamS} and
\eqref{co} imply that $\eta(f) \in U$, and so we can define:

$$ \beta: \{f \in A:\ \s(f)\leq \delta \} \to \Sg, \ \ \beta(f)= P
\circ \eta (f),$$ where $P: U \to \Sg$ is the orthogonal
projection.

Now, let us  check that $\psi(f)=(\beta(f), \sigma(f))$ satisfies
the conditions given by Proposition \ref{covering}. If $\sigma(f)
\leq \delta$, then $\beta(f) \in co [ S(f)] \cap \Sg$. Therefore,
$d(\beta(f), S(f)) < (R_0+1)\sigma(f)$. Take any $p \in S(f)$.
Recall that $R_0 = 3R$ and that $\s(f) \leq 3\s(x,f)<3 \s(f)$ for
any $x \in S(f)$: it is easy to conclude then {\it a)} and {\it
b)}.

If $\sigma(f) \geq \delta$, $\beta$ is not defined. Observe that
{\it a)} is then satisfied for any $\beta \in \Sg$.

\subsection{Proof of Proposition \ref{mt}}

First of all, we will need the following technical lemma:

\begin{lem} \label{technical} There exists $C>0$ such that for any $x \in \Sg$,
$d>0$ small,

$$ \left | \fint_{B_x(d)} u\, dV_g - \fint_{\partial B_x(d)} u\,
dS_g\right | \leq C \left ( \int_{B_{x}(d)} |\nabla u|^2 \, dV_g
\right )^{1/2}.$$

Moreover, given $r\in (0,1)$, there exists $C=C(r, \Sg)>0$ such
that for any $x_1$, $x_2 \in \Sg$, $d>0$ with $B_1 = B_{x_1}(r d)
\subset B_2 = B_{x_2}(d)$, then:

$$ \left | \fint_{B_1} u\, dV_g - \fint_{B_2} u\,
dV_g\right | \leq C \left ( \int_{B_2} |\nabla u|^2 \, dV_g \right
)^{1/2}.$$

\end{lem}

\begin{pf} The existence of such a constant $C$ is given just by the $L^1$
embedding of $H^1$ and trace inequalities. Moreover, $C$ is
independent of $d$ since both inequalities above are dilation
invariant.

\end{pf}

In view of the statement of Proposition \ref{covering}, we now
deduce a Moser-Trudinger type inequality for small balls, and also
for annuli with small internal radius. Those inequalities are in
the core of the proof of Proposition \ref{mt}, and are contained
in the following two lemmas. The first one uses a dilation
argument:

\begin{lem} \label{ball} For any $\e>0$ there exists $C=C(\e)>0$ such that
\begin{eqnarray*}
  (1+\e)\int_{B_p(s)}Q(u_1, u_2) \  dV_g + C & \geq & 4 \pi \left (\log \int_{B_p(s/2)}
  e^{u_1} dV_g +\log \int_{B_p(s/2)} e^{u_2} dV_g \right)  \\
  & - & 4 \pi  (\bar{u}_1(s) + \bar{u}_2(s) + 4 \log s ),
\end{eqnarray*}
for any $u \in H^1(\Sg), \ p \in \Sg$, $s>0$ small and for
$\bar{u}_i(s)= \dis \fint_{B_p(s)} u_i\, dV_g$.
\end{lem}

\begin{pf}
For $s > 0$ smaller than the injectivity radius, the result
follows easily from Proposition \ref{p:MTbd}, for some $C=C(s,
\e)$. Then, we need to prove that the constant $C$ can be taken
independent of $s$ as $s \to 0$. Notice that, as $s \to 0$ we
consider quantities defined on smaller and smaller geodesic balls
$B_q(\varsigma)$ on $\Sg$. Working in normal geodesic coordinates
at $q$, gradients, averages and the volume element will resemble
Euclidean ones. If we assume that near $q$ the metric of $\Sg$ is
flat, we will get negligible error terms which will be omitted for
reasons of brevity.

To prove the lemma, we simply make a dilation of the pair $(u_1, u_2)$ of the form:
$$ v_i(x)= u_i(s x+ p).$$ 
From easy computations there follows:
$$ \int_{B(p,s)} Q(u_1,u_2) \, dV_g= \int_{B(0,1)} Q(v_1, v_2) \, dV_g,  $$
$$ \bar{u}_i(s)=  \fint_{B(0,1)} v_i \, dV_g,  $$
$$ \int_{B_p(s/2)} e^{u_i} \, dV_g =
s^{2} \int_{B(0,1/2)} e^{v_i}\, dV_g.$$ Applying Proposition
\ref{p:MTbd} to the pair $(v_1, v_2)$, we conclude the proof of
the lemma.

\end{pf}

The next lemma gives us an estimate of the quadratic form $Q$ on annuli by
using the Kelvin transform. This transformation is indeed very natural in this
framework, see Remark \ref{r:kelvin} for a more detailed discussion.

\begin{lem} \label{annulus} Given $\e>0$, there exists a fixed $r_0>0$ (depending only on $\Sg$ and $\e$)
satisfying the following property: for any $r \in (0, r_0)$ fixed, there
exists $C=C(r, \e)>0$ such that, for any $(u_1,u_2) \in H^1(\Sg)$
with $u_i=0$ in $\partial B_p(2 r)$,

$$ \int_{A_p(s/2,2 r)}Q(u_1, u_2) \  dV_g + \e \int_{B_p(2 r)}Q(u_1, u_2) \
dV_g + C \geq $$$$ 4 \pi \left (\log \int_{A_p(s, r)} e^{u_1} dV_g
+\log \int_{A_p(s,r)} e^{u_2} dV_g  + (\bar{u}_1(s) + \bar{u}_2(s)
+ 4 \log s )(1+\e)\right ) + C,
$$ with $\ p \in \Sg$, $s\in(0,r)$ and $\bar{u}_i(s)= \fint_{B_p(s)} u_i\,
dV_g$.
\end{lem}

\begin{pf}
As in the proof of Lemma \ref{ball}, we need to show that $C$ is
independent of $s$ as $s \to 0$. By taking $r_0$ small enough,
also here the metric becomes close to the Euclidean one. Reasoning
as in the proof of Lemma \ref{ball}, we can then assume that the
metric is flat around $p$.

We can define the Kelvin transform:
$$ K : A_p(s/2, 2r) \to A_p(s/2, 2r), \qquad  K(x)= p+ r s \frac{x-p}{\ |x-p|^2}.$$

Observe that $K$ maps the interior boundary of $A_p(s/2, 2r)$ onto the exterior one
and viceversa, and fixes the set $\partial B_p(\sqrt{s\, r})$. Let
us define the functions $ \hat{u}_i \in H^1(B_p(2r))$ as:
$$\hat{u}_i(x)= \left \{ \begin{array}{ll} u_i(K(x)) - 4 \log |x-p| & \mbox{ if } |x-p| \geq s/2, \\
-4 \log (s/2)& \mbox{ if } |x-p| \leq s/2. \end{array} \right.$$

Our goal is to apply the Moser-Trudinger inequality given by
Proposition \ref{p:MTbd} to $(\hat{u}_1, \hat{u}_2)$. In order to
do so, let us compute:
\begin{equation} \label{exp} \int_{A_p(s, r)} e^{\hat{u}_i} \, dV_g =  \int_{A_p(s,
r)} e^{u_i(K(x))} |x-p|^{-4} \, dV_g =\frac{1}{s^2 r^2}
\int_{A_p(s, r)} e^{u_i(x)} \, dV_g,\end{equation} since the
Jacobian of $K$ is $J(K(x)) = - r^2 s^2 |x-p|^{-4}$.

Moreover, by Lemma \ref{technical}, we have that:

$$ \left | \fint_{B_p(2r)} \hat{u}_i \, dV_g - \fint_{\partial B_p(r)} \hat{u}_i \,
dS_x \right | \leq C \left (\int_{B_p(2r)} |\n \hat{u}_i| ^2 \,
dV_g  \right )^{1/2} \leq C + \e \int_{B_p(2r)} |\n \hat{u}_i| ^2
\, dV_g .$$

By using again a change of variables,

$$ \fint_{\partial B_p(r)} \hat{u}_i \,
dS_x = \fint_{\partial B_p(s)} u_i \, dS_x - 8 \pi r \log r.$$
Therefore,
\begin{equation} \label{media} \left | \fint_{B_p(2r)} \hat{u}_i \,
dV_g- \bar{u}_i(s) \right | \leq C +  \e \int_{B_p(2r)} |\n
\hat{u}_i| ^2 \, dV_g + \e \int_{B_p(s)} |\n u_i| ^2 \, dV_g.
\end{equation}

Let us now estimate the gradient terms. For $|x-p|\geq s/2$,

$$|\n \hat{u}_i(x)|^2 = |\n u_i(K(x))|^2 \frac{s^2 r^2}{|x-p|^4} +
\frac{16}{|x-p|^2} + 8 \nabla u (K(x)) \cdot \frac{x-p}{|x-p|^4}s
r.$$

Therefore,

$$\int_{B_p(2r)} |\n \hat{u}_i(x)|^2 \, dV_g= \int_{A_p(s/2, 2 r)} |\n \hat{u}_i(x)|^2 \, dV_g =
\int_{A_p(s/2, 2 r)} |\n u_i(K(x))|^2 \frac{s^2 r^2}{\ |x-p|^4} \,
dV_g +$$$$ 16 \int_{A_p(s/2, 2 r)} \frac{dV_g}{|x-p|^2} + 8
\int_{A_p(s/2, 2 r)} \nabla u_i (K(x)) \cdot \frac{x-p}{|x-p|^4}\
s  \, r \, dV_g =
$$$$ \int_{A_p(s/2, 2 r)} |\n u_i(x)|^2 \, dV_g
+ 32 \pi (\log (2r) - \log (s/2)) + 8 \int_{A_p(s/2, 2 r)} \nabla
u_i (K(x)) \cdot \frac{K(x)-p}{|K(x)-p|^2} \ \frac{s^2 r^2}{\
|x-p|^{4}}\, dV_g =
$$
$$ \int_{A_p(s/2, 2 r)} |\n u_i(x)|^2 \, dV_g + 32 \pi (\log (2r) - \log (s/2)) +
8 \int_{A_p(s/2, 2 r)} \nabla u_i (x) \cdot \frac{x-p}{\
|x-p|^{2}} \, dV_g = $$ $$ \int_{A_p(s/2, 2 r)} |\n u_i(x)|^2 \,
dV_g + 32 \pi (\log (2r) - \log (s/2)) - 16 \pi \fint_{\partial
B_p(s/2)} u_i \, dS_x.
$$

In the last equality we have used integration by parts. By using
again Lemma \ref{technical},

\begin{equation} \label{dirich} \left | \int_{B_p(2 r)} |\n \hat{u}_i(x)|^2 \,
dV_g - \int_{A_p(s/2, 2 r)} |\n u_i(x)|^2 \, dV_g + 32 \pi \log s
+ 16 \pi \bar{u}_i(s) \right | \leq C  + \e \int_{B_p(s)} |\n
{u}_i| ^2 \, dV_g.\end{equation}

Regarding the mixed term $\n \hat{u}_1 \cdot \n \hat{u}_2$, we
have that for $|x-p|\geq s/2$,
$$\n \hat{u}_1(x) \cdot \n \hat{u}_2(x) = \n u_1(K(x)) \cdot \n u_2(K(x)) \frac{s^2 r^2}{|x-p|^4} +
\frac{16}{|x-p|^2} + \frac{4sr}{|x-p|^4} (\n u_1(K(x)) + \n
u_2(K(x))\cdot (x-p).$$ Reasoning as above, we obtain the
estimate:

\begin{equation} \label{dirich2}  \begin{array}{c} \left | \displaystyle \int_{B_p(2 r)} \n \hat{u}_1(x) \cdot \n \hat{u}_2(x) \,
dV_g - \int_{A_p(s/2, 2 r)} \n u_1(x) \cdot \n u_2(x) \, dV_g + 32
\pi
\log s + 8 \pi \bar{u}_1(s) + 8 \pi \bar{u}_2(s) \right | \leq \\
\\  C  + \e \displaystyle \int_{B_p(s)} \left( |\n {u}_1| ^2 + |\n
{u}_2|^2 \right )\, dV_g. \end{array}\end{equation}

We now apply Proposition \ref{p:MTbd} to $(\hat{u}_1, \hat{u}_2)$
and use the estimates \eqref{exp}, \eqref{media}, \eqref{dirich}
and \eqref{dirich2}, to obtain:

$$ 4 \pi \left [ \log \int_{A_p(s,r)} e^{u_1} \, dV_g+
 \log \int_{A_p(s,r)} e^{u_2} \, dV_g -(4 \log s + \bar{u}_1(s) + \bar{u}_2(s)) \right ]   \leq $$$$
 4 \pi \left [ \log \int_{A_p(s,r)}
e^{\hat{u}_1} \, dV_g +
 \log \int_{A_p(s,r)}  e^{\hat{u}_2}\, dV_g - \fint_{B_p(2r)} (\hat{u}_1 + \hat{u}_2) \right ]
 +$$$$
 \e \int_{B_p(2r)} \left(|\nabla \hat{u}_1|^2 +
|\nabla \hat{u}_2|^2 \right) \, dV_g + \e \int_{B_p(s)} \left (
|\n u_1|^2 + |\n u_2|^2 \right )\, dV_g + C \leq
$$
$$ (1+C\e) \int_{B_p(2r)} Q(\hat{u}_1,\hat{u}_2)\, dV_g + \e \int_{B_p(s)} \left (
|\n u_1|^2 + |\n u_2|^2 \right )\, dV_g +C \leq $$
$$ (1+C\e) \left [ \int_{A_p(s/2,2r)} Q(u_1,u_2)\, dV_g - 8 \pi ( 4
\log s + \bar{u}_1(s) + \bar{u}_2(s) ) \right ] + \e \int_{B_p(s)}
\left( |\nabla u_1|^2 + |\nabla u_2|^2 \right) \, dV_g+C.$$

%

By renaming $\e$ conveniently, we conclude the proof.

\end{pf}

\begin{rem} The term $\bar{u}(s) + 2
\log s$ has an easy interpretation; by the Jensen inequality we
have the estimate
$$ \log \int_{B_p(s)} e^u \, dV_g = \log \left (|B_p(s)| \fint_{B_p(s)} e^u \, dV_g \right) \geq \bar{u}(s) + 2
\log s - C.$$

\end{rem}

\begin{rem}\label{r:kelvin} The transformation $K$ is used to exploit
the geometric properties of the problem, in order to gain as much
control as possible on the exponential terms. From the formulas in
\cite{jw2} one has that both components of the entire solutions of
the Toda system in $\R^2$ decay at infinity at the rate $- 4 \log |x|$. In
this way, the Kelvin transform brings these functions to (nearly)
constants at the origin, giving a sort of optimization in the
Dirichlet part. The minimal value of Dirichlet energy to obtain
concentration of volume at a scale $s$ (as in the statement of Lemma
\ref{annulus}) is then transformed into a boundary integral which cancels
exactly the extra terms in Lemma \ref{ball} due to the $s$-dilation.
\end{rem}

\begin{rem}

Lemmas \ref{ball} and \ref{annulus}, together with Proposition
\ref{covering}, give a precise idea of the proof. Indeed, assume
that for some $p\in \Sg$, $\s >0$:

\begin{equation} \label{1} \int_{B_{p}(\s)} e^{u_i} dV_g \geq \tau \int_\Sg e^{u_i}
dV_g,
  \ i=1,2;\end{equation}
\begin{equation} \label{2} \int_{B_{p}(R \s)^c} e^{u_i} dV_g \geq \tau \int_\Sg e^{u_i}
dV_g, \ i=1,2.
\end{equation}
If we sum the inequalities given by Lemmas \ref{ball} and
\ref{annulus}, the term $\bar{u}_1(\s) + \bar{u}_2(\s) + 4 \log
\s$ cancels and we deduce the estimate of Proposition \ref{mt}.

The problem is that when $\psi \left( \frac{e^{u_1}}{\int_{\Sigma}
e^{u_1} dV_g} \right )= \psi \left( \frac{e^{u_2}}{\int_{\Sigma}
e^{u_2} dV_g} \right )$ we do not really have \eqref{1}, \eqref{2}
around the same point $p$. Moreover, $u_i$ needs not be zero on
the boundary of a ball, as requested in Proposition \ref{annulus}.
Some technical work is needed to deal with those difficulties.

\end{rem}

We now prove Proposition \ref{mt}. Fixed $\e>0$, take $R>1$
(depending only on $\e$) and let $\psi$ be the continuous map
given by Proposition \ref{covering}. Fix also $\delta>0$ small
(which will depend only on $\e$, too).

Let $u_1$ and $u_2$ be two functions in $H^1(\Sg)$ with
$\int_{\Sg} u_{i} \, dV_g =0$, such that:

$$\psi \left( \frac{e^{u_1}}{\int_{\Sigma} e^{u_1} dV_g} \right
)= \psi \left( \frac{e^{u_2}}{\int_{\Sigma} e^{u_2} dV_g} \right
)= (\beta, \sigma) \in \ov{\Sg}_{\d}. $$

If $\sigma \geq \frac{\delta}{R^2}$, then Proposition
\ref{p:imprc} yields the result. Therefore, assume $\sigma <
\frac{\delta}{R^2}$; Proposition \ref{covering} implies the
existence of $\tau>0$, $p_1,\ p_2 \in \Sg$ satisfying:
\begin{equation} \label{dentro} \int_{B_{p_i}(\s)} e^{u_i} dV_g \geq \tau \int_\Sg
e^{u_i} dV_g,
  \ i=1,2;\end{equation}
\begin{equation}
  \int_{B_{p_i}(R \s)^c} e^{u_i} dV_g \geq \tau \int_\Sg e^{u_i} dV_g \ i=1,2;
\end{equation}
$$d(p_1, p_2)
  \leq (6R+2) \s;
$$

The proof will be divided into two cases:

\bigskip

\noindent {\bf CASE 1:} Assume that:

\begin{equation} \label{fuera} \int_{A_{p_i}(R\s, \delta)} e^{u_i} dV_g \geq \tau/2 \int_{\Sg} e^{u_i}
dV_g.\end{equation}

In order to be able to apply Lemma \ref{annulus}, we need to
modify our functions outside a certain ball. Choose $k \in \N$, $k
\leq 2 \e^{-1}$, such that:

$$ \int_{A_{p_1}(2^{k-1} \delta, 2^{k+1} \delta)} \left ( |\n u_1|^2 + |\n u_2|^2 \right )\, dV_g
\leq \e \int_{\Sg} \left (|\n u_1|^2 + |\n u_2|^2 \right )\, dV_g.
$$

We define $\tilde{u}_i \in H^1(\Sg)$ by:

$$  \left \{ \begin{array}{ll} \tilde{u}_i(x) = u_i(x) & x \in
B_{p_1}(2^k \delta), \\ \Delta \tilde{u}_i(x) =0 & x \in
A_{p_1}(2^k \delta, 2^{k+1} \delta), \\ \tilde{u}_i(x) = 0 & x \notin
B_{p_1}(2^{k+1} \delta).
\end{array} \right.  $$

Since we plan to apply Lemma \ref{annulus} to $(\tilde{u}_1,
\tilde{u}_2)$, we need to choose $\delta$ small enough so that
$2^{3\e^{-1}}\delta < r_0$, where $r_0$ is given by that Lemma.

It is easy to check, by using Lemma \ref{technical}, that

$$ \int_{A_{p_1}(2^{k} \delta, 2^{k+1} \delta)} \left (|\n \tilde{u}_1|^2 + |\n \tilde{u}_2|^2 \right ) \, dV_g \leq $$$$
C \int_{A_{p_1}(2^{k-1} \delta, 2^{k} \delta)} \left (|\n u_1|^2 +
|\n u_2|^2 \right ) \, dV_g \leq C \e \int_{\Sg} \left (|\n u_1|^2
+ |\n u_2|^2 \right )\, dV_g,  $$ where $C$ is a universal
constant.

\bigskip

\noindent {\bf Case 1.1:} $d(p_1, p_2) \leq R^{\frac 12} \s$.

\medskip
By applying Lemma \ref{ball} to $u_i$ for $p=p_1$ and $s= 2
(R^{1/2}+1)\sigma$, and taking into account \eqref{dentro}, we
obtain:

\begin{equation} \label{dentro11} \begin{array}{c}(1+\e) \dis \int_{B_p(s)}Q(u_1, u_2) \  dV_g + C \geq  \\ \\4 \pi \left
(\log \dis \int_{B_p(s/2)} e^{u_1} dV_g +\log \int_{B_p(s/2)}
e^{u_2} dV_g  - (\bar{u}_1(\s) + \bar{u}_2(\s) + 4 \log \s )
\right ) \geq
\\ \\ 4 \pi \left (\log \dis \int_{\Sg} e^{u_1} dV_g +\log \int_{\Sg}
e^{u_2} dV_g  - (\bar{u}_1(\s) + \bar{u}_2(\s) + 4 \log \s )-C
\right ),\end{array}
\end{equation} where $\bar{u}_i(\s)= \fint_{B_p(\s)} u_i\, dV_g$.

We now apply Lemma \ref{annulus} to $\tilde{u}_i$ for $p=p_1$,
$s'=4(R^{1/2}+1)\s$ and $r = 2^{k+1}\delta$:

\begin{equation} \label{fuera11} \begin{array}{c}
\dis \int_{A_p(s'/2,2 r)}Q(\tilde{u}_1, \tilde{u}_2) \  dV_g  + \e
\int_{\Sg} \left( |\n u_1|^2 + |\n u_2|^2 \right )\, dV_g+ C \geq \\ \\
4 \pi \left (\log \dis \int_{A_p(s', r)} e^{\tilde{u}_1} dV_g
+\log \int_{A_p(s',r)} e^{\tilde{u}_2} dV_g  + (\bar{u}_1(\s) +
\bar{u}_2(\s) + 4 \log \s )(1+\e)\right ).
\end{array} \end{equation}

Taking into account \eqref{fuera}, we conclude:

\begin{equation} \label{fuera11-bis} \begin{array}{c}
\dis \int_{A_p(s'/2,2 r)}Q(\tilde{u}_1, \tilde{u}_2) \  dV_g  + \e
\int_{\Sg} \left( |\n u_1|^2 + |\n u_2|^2 \right ) \, dV_g+ C \geq \\
\\4 \pi \left (\log \dis \int_{\Sg} e^{u_1} dV_g +\log \int_{\Sg}
e^{{u}_2} dV_g  + (\bar{u}_1(\s) + \bar{u}_2(\s) + 4 \log \s
)(1+\e)\right ).\end{array} \end{equation}

Combining \eqref{dentro11} and \eqref{fuera11-bis} we obtain our
result (after properly renaming $\e$).

\bigskip

\noindent {\bf Case 1.2:} $d(p_1, p_2) \geq R^{\frac 12} \s$
and $\displaystyle \int_{B_{p_1}(R^{\frac 13} \s)} e^{u_2} dV_g
\geq \tau/4 \int_\Sg e^{u_2} dV_g$.

\medskip

Here we argue as in Case 1.1: as a first step, we apply Lemma
\ref{ball} to $(u_1, u_2)$ for $p=p_1$ and $s = 2 (R^{1/3}+1)\s$.
Then, we use Lemma \ref{annulus} with $(\tilde{u}_1, \tilde{u}_2)$
for $p=p_1$, $s'= 4 (R^{1/3}+1)\s$ and $r=2^{k+1} \delta$.

\bigskip

\noindent {\bf Case 1.3:} $d(p_1, p_2) \geq R^{\frac 12} \s$
and $\displaystyle \int_{B_{p_2}(R^{\frac 13} \s)} e^{u_1} dV_g
\geq \tau/4 \int_\Sg e^{u_1} dV_g$.

\medskip

This case can be treated as in Case 1.2, by just interchanging the
indices $1$ and $2$.

\bigskip

\noindent {\bf Case 1.4:} $d(p_1, p_2) \geq R^{\frac 12} \s$,
$\displaystyle \int_{B_{p_2}(R^{\frac 13} \s)} e^{u_1} dV_g \leq
\tau/4 \int_\Sg e^{u_1} dV_g$ and $\displaystyle
\int_{B_{p_1}(R^{\frac 13} \s)} e^{u_2} dV_g \leq \tau/4 \int_\Sg
e^{u_2} dV_g$.

\medskip

Here we need to use again some harmonic lifting of our functions.
Take $n \in \N$, $n \leq 2 \e^{-1}$ so that

$$ \sum_{i=1}^2 \int_{A_{p_i}(2^{n-1} \s, 2^{n+1} \s )} \left( |\n u_1|^2 + |\n u_2|^2 \right ) \, dV_g
\leq \e \int_{\Sg} \left (|\n u_1|^2 + |\n u_2|^2 \right )\, dV_g,
$$ where we have chosen $R$ so that $2^{3\e^{-1}} <R^{1/3}$.

We define the function $v_i$ of class $H^1$ by:
$$  \left \{ \begin{array}{ll} v_i(x) = u_i(x) & x \in
B_{p_1}(2^{n} \s) \cup B_{p_2}(2^{n} \s), \\
\Delta v_i(x) =0 & x \in A_{p_1}(2^{n} \s, 2^{n+1} \s) \cup A_{p_2}(2^{n}\s, 2^{n+1} \s), \\
v_i(x) = \bar{u}_i(\s) & x \notin B_{p_1}(2^{n+1} \s) \cup B_{p_2}(2^{n+1}
\s).
\end{array} \right.  $$

Again,

$$ \sum_{i=1}^2 \int_{A_{p_i}(2^{n} \s, 2^{n+1} \s)} \left( |\n v_1|^2 + |\n v_2|^2 \right) \, dV_g \leq $$$$
C \sum_{i=1}^2 \int_{A_{p_i}(2^{n-1}  \s, 2^{n+1} \s)} \left (|\n
u_1|^2 + |\n u_2|^2 \right ) \, dV_g \leq C \e \int_{\Sg} \left
(|\n u_1|^2 + |\n u_2|^2 \right )\, dV_g,  $$ where $C$ is a
universal constant.

We now apply Lemma \ref{ball} to $(v_1,v_2)$ with $p=p_1$ and
$s=2(6R+2)\s$, and take into account \eqref{dentro}:

\begin{equation} \label{dentro14} \begin{array}{c}
\dis \int_{ B_{p_1}(2^{n} \s) \cup B_{p_2}(2^{n} \s)} Q(u_1,u_2) \, dV_g
+ C \e \int_{\Sg} \left ( |\n u_1|^2 + |\n u_2|^2 \right) \, dV_g +C \geq \\ \\
 (1+\e)\dis \int_{B_p(s)}Q(v_1, v_2) \  dV_g +C \geq \\ \\ 4 \pi \left (\log
\dis \int_{B_p(s/2)} e^{v_1} dV_g +\log \int_{B_p(s/2)} e^{v_2}
dV_g - (\bar{u}_1(s) + \bar{u}_2(s) + 4 \log s ) \right ) \geq \\
\\ 4 \pi \left (\log \dis \int_{\Sg} e^{u_1} dV_g +\log \int_{\Sg}
e^{u_2} dV_g - (\bar{u}_1(s) + \bar{u}_2(s) + 4 \log s ) \right
)-C.
\end{array} \end{equation}

Now, we define $w_i \in H^1(\Sg)$ as:

$$  \left \{ \begin{array}{ll} w_i(x) = \bar{u}_i(\s) & x \in
B_{p_1}(2^{n} \s) \cup B_{p_2}(2^{n} \s), \\
\Delta w_i(x) =0 & x \in A_{p_1}(2^{n} \s, 2^{n+1} \s) \cup A_{p_2}(2^{n}\s, 2^{n+1} \s), \\
w_i(x) = \tilde{u}_i & x \notin B_{p_1}(2^{n+1} \s) \cup B_{p_2}(2^{n+1}
\s).
\end{array} \right.  $$

As before,

$$ \sum_{i=1}^2 \int_{A_{p_i}(2^{n} \s, 2^{n+1} \s)} \left (|\n w_1|^2 + |\n w_2|^2 \right ) \, dV_g \leq $$$$
C \sum_{i=1}^2 \int_{A_{p_i}(2^{n-1}  \s, 2^{n+1} \s)} \left (|\n
u_1|^2 + |\n u_2|^2 \right ) \, dV_g \leq C \e \int_{\Sg} \left
(|\n u_1|^2 + |\n u_2|^2 \right ) \, dV_g,  $$ where also here $C$
is a universal constant.

We apply Lemma \ref{annulus} to $(w_1,w_2)$ for any point $p'$
such that $d(p', p_1) = \frac 1 2 R^{1/3}\s$, $s'= \s$ and $r =
2^{k+1} \delta$:
\begin{equation} \label{fuera14} \begin{array}{c}\dis
\int_{ ( B_{p_1}(2^{n+1} \s) \cup B_{p_2}(2^{n+1} \s))^c}
Q(u_1,u_2) \, dV_g + C \e \int_{\Sg} \left (|\n u_1|^2 + |\n
u_2|^2 \right )\, dV_g  +C \geq
\\ \\ (1+\e)\dis \int_{A_{p'}(s'/2,2 r)}Q(w_1, w_2) \  dV_g +C \geq
\\ \\ 4 \pi \left (\log \dis \int_{A_{p'}(s', r)} e^{w_1} dV_g +\log
 \int_{A_{p'}(s',r)} e^{w_2} dV_g  + (\bar{u}_1(\s) +
\bar{u}_2(\s) + 4 \log \s )(1+\e)\right ). \end{array}
\end{equation}

Taking into account \eqref{fuera} and the hypothesis of Case 1.4,

\begin{equation} \label{fuera14-bis} \begin{array}{c}\dis
\int_{ ( B_{p_1}(2^n \s) \cup B_{p_2}(2^n \s))^c} Q(u_1,u_2) \,
dV_g + C \e \int_{\Sg} \left ( |\n u_1|^2 + |\n u_2|^2 \right )\,
dV_g +C \geq
\\ \\ 4 \pi \left (\log \dis \int_{\Sg} e^{u_1} dV_g +\log
 \int_{\Sg} e^{u_2} dV_g  + (\bar{u}_1(\s) +
\bar{u}_2(\s) + 4 \log \s )(1+\e)\right ).
\end{array} \end{equation}

Combining inequalities \eqref{dentro14} and \eqref{fuera14-bis},
we obtain our result.

\bigskip

\noindent {\bf CASE 2:} Assume that for some $i=1,2$:

$$ \int_{B_{p_i}(\delta)^c} e^{u_i} dV_g \geq \tau/2   \int_{\Sg} e^{u_i} dV_g.$$

Without loss of generality, let us consider $i=1$.

Take $\delta'= \frac{\delta}{2^{3/\e}}$. If moreover:
$$ \int_{B_{p_2}(\delta')^c} e^{u_2} dV_g \geq \tau/2,$$
then Proposition \ref{p:imprc} implies the desired inequality. So,
we can assume that: \begin{equation} \label{hia}
\int_{A_{p_2}(R\s, \delta')} e^{u_2} dV_g \geq
\tau/2.\end{equation}

We now apply the whole procedure of Case 1 to $u_1$, $u_2$,
replacing $\delta$ with $\delta'$.

For instance, as in Case 1.1, we would get \eqref{dentro11} and
\eqref{fuera11}. However, here \eqref{fuera11-bis} does not follow
immediately since now we do not know whether:

$$ \int_{A_p(s',r)} e^{u_1} dV_g \geq \alpha \int_{\Sg} e^{u_1} dV_g,$$
for some fixed $\alpha>0$. This is needed to estimate:
$$ \log \int_{A_p(s', r)} e^{\tilde{u}_1} dV_g \geq \log \int_{\Sg} e^{u_1}
dV_g- C,$$ which allows us to obtain \eqref{fuera11-bis}.

By applying the Jensen inequality and Lemma \ref{technical}, we
get:

$$ \log \int_{A_p(s', r)} e^{\tilde{u}_1} dV_g \geq \log \int_{A_p(r/8, r/4)} e^{u_1}
dV_g\geq $$
$$ \log \fint_{A_p(r/8, r/4)} e^{u_1}
dV_g-C \geq \fint_{A_p(r/8, r/4)} u_1 dV_g-C \geq -\e \int_{\Sg}
|\n u_1|^2\, dV_g - C.
$$

Therefore, from \eqref{hia} and \eqref{fuera11} we get:

\begin{equation} \label{fuera2} \begin{array}{c}
\dis \int_{A_p(s'/2,2 r)}Q(\tilde{u}_1, \tilde{u}_2) \  dV_g  +
C\e \int_{\Sg} \left (|\n u_1|^2 + |\n u_2|^2 \right )\,  dV_g+ C \geq \\
\\4 \pi \left (\log \dis \int_{\Sg} e^{{u}_2} dV_g  +
(\bar{u}_1(\s) + \bar{u}_2(\s) + 4 \log \s )(1+\e)\right
).\end{array}
\end{equation}

Now, we apply Proposition \ref{p:MTbd}, to find:
$$
(1+C\e) \int_{B_{p_1}(\delta/2)^c} Q(u_1, u_2) dV_g  +C \geq 4 \pi
\left ( \log \int_{B_{p_1}(\delta)^c}  e^{u_1} dV_g + \log
\int_{B_{p_1}(\delta)^c} e^{u_2} dV_g \right ).
$$
Again here we can use Jensen inequality and the hypothesis of Case
2 to deduce:

\begin{equation} \label{fuera2-bis} \begin{array}{c}
\dis \int_{B_{p_1}(\delta/2)^c} Q(u_1, u_2) dV_g + C\e \int_{\Sg}
\left (|\n u_1|^2 + |\n u_2|^2 \right )\,  dV_g+ C\geq \\ \\4 \pi
\left ( \log \dis \int_{\Sg}  e^{u_1} dV_g \right ).
\end{array}\end{equation}

We conclude now by combining \eqref{fuera2-bis}, \eqref{fuera2}
and \eqref{dentro11}.

We can argue in the same way if we are under the conditions of
Cases 1.2, 1.3 or 1.4.

\begin{rem}\label{puffff} The improved inequality in Proposition \ref{mt} is
consistent with the asymptotic analysis in \cite{jlw}. Here the authors prove that
when both $u_1, u_2$ blow-up at the same rate at the same point, then the
corresponding quantization of conformal volume is $(8\pi, 8\pi)$. On the other
hand when the blow-up rates are different, but occur at the same point, then
the quantization values are $(4\pi, 8\pi)$ or $(8\pi,4\pi)$.
\end{rem}

\section{Min-max scheme}\label{s:4}

\noindent Let $\ov{\Sg}_{\d}$ be as in \eqref{cono}, and let us
set
\begin{equation}\label{eq:DX}
    \ov{D}_{\d} = diag(\ov{\Sg}_{\d} \times \ov{\Sg}_{\d}) := \left\{ (\vartheta_1, \vartheta_2)
    \in \ov{\Sg}_{\d} \times \ov{\Sg}_{\d} \; : \; \vartheta_1 = \vartheta_2 \right\};
    \qquad \quad X = \ov{\Sg}_{\d} \times \ov{\Sg}_{\d} \setminus \ov{D}_{\d}.
\end{equation}
Let $\e > 0$ be such that $\rho_i + \e < 8 \pi$ for $i = 1, 2$, and let $R, \d, \psi$
be as in Proposition \ref{covering}. Consider then the map $\Psi : H^1(\Sg)
\times H^1(\Sg)$ defined in the following way
\begin{equation}\label{eq:Psi}
    \Psi(u_1, u_2) = \left( \psi \left( \frac{e^{u_1}}{\int_{\Sg} e^{u_1} dV_g}
    \right), \psi \left( \frac{e^{u_2}}{\int_{\Sg} e^{u_2} dV_g} \right) \right).
\end{equation}
By Proposition \ref{mt}, and since $C \geq h_i(x) \geq
\frac{1}{C}>0$ for any $x \in \Sg$, we have that $I_\rho(u_1,
u_2)$ is bounded from below for any $(u_1,\ u_2)$ such that
$\Psi(u_1, u_2) \in \ov{D}_{\d}$. Therefore, there exists a large
$L > 0$ such that
\begin{equation}\label{eq:psilow}
    J_\rho(u_1, u_2) \leq - L \quad \Rightarrow \quad \Psi(u_1, u_2) \in X.
\end{equation}

\

\noindent By our definition of $\ov{\Sg}_{\d}$, the set $X$ is not compact: however it
retracts to some compact subset $\mathcal{X}_\nu$, as it is shown in the next result.

\begin{lem}\label{l:retr} For $\nu \ll \d$, define
$$
  \mathcal{X}_{\nu,1} = \left\{ \left( (x_1, t_1), (x_2, t_2) \right) \in X
    \; : \; \left| t_1 - t_2 \right|^2 + d(x_1, x_2)^2 \geq \d^4, \max\{t_1, t_2\}
    < \d, \min\{t_1, t_2\} \in \left[ \nu^2, \nu \right] \right\};
$$
$$
  \mathcal{X}_{\nu,2} = \left\{ \left( (x_1, t_1), (x_2, t_2) \right) \in X
    \; : \; \max\{t_1, t_2\} = \d, \min\{t_1, t_2\}
    \in \left[ \nu^2, \nu \right] \right\},
$$
and set
\begin{equation}\label{eq:retrx}
    \mathcal{X}_{\nu} = \left( \mathcal{X}_{\nu,1}
    \cup \mathcal{X}_{\nu,2} \right) \subseteq X.
\end{equation}
Then there is a retraction $R_{\nu}$ of $X$ onto $\mathcal{X}_{\nu}$.
\end{lem}

\begin{pf} We proceed in two steps. First, we define a deformation of $X$ in itself
satisfying that:

\begin{enumerate} \label{caracola}

\item[a)] either $\max\{t_1, t_2\} < \d$ and $\left| t_1 - t_2
\right|^2 + d(x_1, x_2)^2 \geq \d^4$,

\item[b)] or $\max\{t_1, t_2\} = \d$.

\end{enumerate}

Then another deformation will provide us with the condition
$\min\{t_1, t_2\} \in \left[ \nu^2, \nu \right]$.

Let us consider the following ODE in $(\Sg \times (0,\delta])^2$:
$$
  \frac{d}{ds} \left(
                 \begin{array}{c}
                   x_1(s) \\
                   t_1(s) \\
                   x_2(s) \\
                   t_2(s) \\
                 \end{array}
               \right) =
\left(
  \begin{array}{c}
    (\delta - \max_i \{t_i(s)\}) \n_{x_1} d(x_1(s), x_2(s))^2  \\
    (t_1(s) - t_2(s)) t_1(s) (\delta - t_1(s))  \\
    (\delta - \max_i \{t_i(s)\}) \n_{x_2} d(x_1(s), x_2(s))^2 \\
    (t_2(s) - t_1(s)) t_2(s) (\delta - t_2(s)) \\
  \end{array}
\right).
$$

Notice that if  $\left| t_1 - t_2 \right|^2 + d(x_1, x_2) < \d^4$
(and $\max\{t_1, t_2\} < \d$) then $d(x_1, x_2)$ is small so
$d(x_1, x_2)^2$ is a smooth function on $(\Sg \times \R)^2$, and
the above vector field is well defined. For each initial datum
$(\vartheta_1, \vartheta_2) \in X$ we define $s_{\vartheta_1,
\vartheta_2}\geq 0$ as the smallest value of $s$ for which the
above flow satisfies either a) or b).

To define the first homotopy $H_1(s,\cdot)$ then one can use the
above flow, rescaling in the evolution variable (depending on the
initial datum) as $s \mapsto \tilde{s} = s_{\vartheta_1,
\vartheta_2}s$.

To define the second homotopy, we introduce two cutoff functions
$\chi_1, \chi_2$:
$$ \begin{array}{ll}
  \left\{
      \begin{array}{ll}
        \chi_1(t) = 1 & \hbox{ for } t \leq \nu^2. \\
        \chi \mbox{ is non increasing, } & \\
        \chi_1(t) = -1& \hbox{ for } t \geq \nu,
      \end{array}
    \right.
    &
  \left\{
      \begin{array}{ll}
        \chi_2(t) = 1 & \hbox{ for } t \leq \delta/2, \\
        \chi_2(t)= 2 \left( 1-\frac{t}{\delta}\right) & t \in (\delta/2, \delta),  \\
        \chi_2(t) = 0& \hbox{ for } t \geq \delta,
      \end{array}
    \right.
    \end{array}
$$
and consider the following ODE
$$
  \frac{d}{ds} \left(
                 \begin{array}{c}
                   t_1(s) \\
                   t_2(s) \\
                 \end{array}
               \right) = \left(
                           \begin{array}{c}
                             \chi_1(\min_i\{t_i(s)\})
                             \chi_2(t_1(s)) \\
                             \chi_1(\min_i\{t_i(s)\})
                             \chi_2(t_2(s))
                           \end{array}
                         \right).
$$
As in the previous case, there exists $\hat{s}_{\vartheta_1,
\vartheta_2}$ such that the condition $\min_i t_i \in [\nu^2,
\nu]$ is reached for $s = \hat{s}_{\vartheta_1, \vartheta_2}$, and
one can define the homotopy $H_2$ rescaling in $s$
correspondingly. Observe that along the homotopy $H_2$ the
distance $| t_1 - t_2 |$ is non decreasing if $|t_1-t_2| \leq
\d/4$.

The concatenation of the homotopies $H_1$ and $H_2$ gives the
desired conclusion. Note that both $H_1$ and $H_2$, by the way
they are constructed, preserve the quotient relations in the
definition of $X$.
\end{pf}

\

\noindent We next construct a family of test functions parameterized by $\mathcal{X}_{\nu}$
on which $J_{\rho}$ attains large negative values. For $(\vartheta_1, \vartheta_2) = \left( (x_1, t_1), (x_2, t_2) \right) \in \mathcal{X}_\nu$ define
\begin{equation}\label{eq:test}
    \var_{(\vartheta_1, \vartheta_2)}(y) = \left( \var_1(y), \var_2(y) \right),
\end{equation}
where we have set
\begin{equation}\label{eq:varvar12}
    \var_1(y) = \log  \frac{1 + \tilde{t}_2^2
    d(x_2,y)^2}{\left( 1 + \tilde{t}_1^2 d(x_1,y)^2 \right)^2}, \qquad \qquad
    \var_2(y) =  \log \frac{1 +
    \tilde{t}_1^2 d(x_1,y)^2}{\left( 1 + \tilde{t}_2^2 d(x_2,y)^2 \right)^2},
\end{equation}
with
\begin{equation}\label{eq:tti}
    \tilde{t}_1 = \tilde{t}_1(t_1) = \left\{
                    \begin{array}{ll}
                      \frac{1}{t_1}, & \hbox{ for } t_1 \leq \frac{\d}{2}, \\
                      - \frac{4}{\d^2} (t_1-\d) & \hbox{ for } t_1 \geq \frac{\d}{2};
                    \end{array}
                  \right. \qquad
\tilde{t}_2 = \tilde{t}_2(t_2) = \left\{
                    \begin{array}{ll}
                      \frac{1}{t_2}, & \hbox{ for } t_2 \leq \frac{\d}{2}, \\
                      - \frac{4}{\d^2} (t_2-\d) & \hbox{ for } t_2 \geq \frac{\d}{2}.
                    \end{array}
                  \right.
\end{equation}
Notice that, by our choices of $\tilde{t}_1, \tilde{t}_2$, this
map is well defined on $\mathcal{X}_\nu$ (especially for what
concerns the identifications in $\ov{\Sg}_\d$). We have then the
following result.

\begin{lem}\label{l:integrals} For $\nu$ sufficiently small, and for $(\vartheta_1,
\vartheta_2) \in \mathcal{X}_{\nu}$, there exists a constant
$C=C(\d,\Sg) > 0$, depending only on $\Sg$ and $\d$, such that
\begin{equation}\label{eq:inttot}
    \frac{1}{C} \frac{t_i^2}{t_j^2} \leq \int_\Sg e^{\var_i(y)}
    dV_g(y) \leq C \frac{t_i^2}{t_j^2}, \qquad \quad i \neq j;
\end{equation}
\end{lem}

\begin{pf} First, we notice that by an elementary change of variables
\begin{equation}\label{eq:intfalt}
    \int_{\R^2} \frac{1}{\left( 1 + \l^2 |x|^2 \right)^2} dx = \frac{C_0}{\l^2};
     \qquad \quad \l > 0
\end{equation}
for some fixed positive constant $C_0$. We distinguish next the two cases
\begin{equation}\label{eq:alt1}
    |t_1 - t_2| \geq \d^3 \qquad \quad \hbox{ and } \qquad \quad |t_1 - t_2| < \d^3.
\end{equation}
In the first alternative, by the definition of $\mathcal{X}_{\nu}$ and by the
fact that $\nu \ll \d$, one of the $t_i$'s belongs to $[\nu^2, \nu]$,
while the other is greater or equal to $\frac{\d^3}{2}$.

If $t_1 \in [\nu^2,
\nu]$ and if $t_2 \geq \frac{\d^3}{2}$ then the function $1 + \tilde{t}_2^2
d(x_2,y)^2$ is bounded above and below by two positive constants depending only
on $\Sg$ and $\d$. Therefore, working in geodesic normal coordinates centered
at $x_1$ and using \eqref{eq:intfalt} we obtain
$$
    \frac{t_1^2}{C} \leq \frac{1}{C \tilde{t}_1^2} \leq \int_\Sg
    e^{\var_1(y)} dV_g(y) \leq \frac{C}{\tilde{t}_1^2} \leq C t_1^2.
$$
If instead $t_2 \in [\nu^2, \nu]$ and if $t_1 \geq \frac{\d^3}{2}$
then the function $1 + \tilde{t}_1^2 d(x_1,y)^2$ is bounded above and below by
two positive constants depending only on $\Sg$ and $\d$, hence one finds
$$
  \int_\Sg e^{\var_1(y)} dV_g(y) \geq \frac{1}{C}
   \int_\Sg (1 + \tilde{t}_2^2 d(x_2,y)^2) dV_g(y) \geq
 \frac{\tilde{t}_2^2}{C} = \frac{1}{C t_2^2},
$$
and similarly
$$
  \int_\Sg e^{\var_1(y)} dV_g(y) \leq C
   \int_\Sg (1 + \tilde{t}_2^2 d(x_2,y)^2) dV_g(y) \leq
 C \tilde{t}_2^2 = \frac{C}{t_2^2}.
$$
In both the last two cases we then obtain the conclusion.

Suppose now that $|t_1 - t_2| < \d^3$: then by the definition of $\mathcal{X}_{\nu}$
we have that $d(x_1, x_2) \geq \frac{\d^2}{2}$ and that $t_1, t_2 \leq \nu
+ \d^3$.  Then, from \eqref{eq:intfalt} and some elementary estimates we derive
$$
  \int_\Sg e^{\var_1(y)} dV_g(y) \geq \int_{B_{x_1}(\delta^3)} e^{\var_1(y)} dV_g(y)
  \geq \frac{1}{C} \frac{1 + \tilde{t}_2^2
 d(x_1,x_2)^2}{\tilde{t}_1^2} \geq \frac{1}{C} \frac{t_1^2}{t_2^2}.
$$
By the same argument we obtain
$$
  \int_{B_{x_1}(\delta^3)} e^{\var_1(y)} dV_g(y) \leq C \frac{1 +
  \tilde{t}_2^2 d(x_1,x_2)^2}{\tilde{t}_1^2} \leq C \frac{t_1^2}{t_2^2}.
$$
Moreover, we have
$$
  \int_{(B_{x_1}(\delta^3))^c} e^{\var_1(y)} dV_g(y) \leq \frac{C}{\tilde{t}_1^4}
   \int_{(B_{x_1}(\delta^3))^c} (1 + \tilde{t}_2^2 d(x_2,y)^2) dV_g(y)
   \leq C \frac{t_1^4}{t_2^2}.
$$
This concludes the proof.
\end{pf}

\begin{lem}\label{l:dsmallIlow} For $(\vartheta_1, \vartheta_2) \in
\mathcal{X}_{\nu}$, let $\var_{(\vartheta_1, \vartheta_2)}$ be
defined as in the above formula. Then
$$
  J_{\rho}(\var_{(\vartheta_1, \vartheta_2)}) \to - \infty \quad
  \hbox{ as } \nu \to 0 \qquad \quad \hbox{ uniformly for }
   (\vartheta_1, \vartheta_2) \in \mathcal{X}_{\nu}.
$$
\end{lem}

\begin{pf} The statement follows from Lemma \ref{l:integrals} once the
following three estimates are shown
\begin{equation}\label{eq:estQ}
    \int_{\Sg} Q\left( \var_{(\vartheta_1, \vartheta_2)} \right) dV_g \leq
    8 \pi (1 + o_{\d}(1)) \log \frac{1}{t_1} +
    8 \pi (1 + o_{\d}(1)) \log \frac{1}{t_2};
\end{equation}
\begin{equation}\label{eq:estlin1}
    \fint_{\Sg} \var_1 dV_g = 4 (1 + o_{{\d}}(1))
   \log t_1 - 2 (1 + o_{{\d}}(1)) \log t_2;
\end{equation}
\begin{equation}\label{eq:estlin2}
    \fint_{\Sg} \var_2 dV_g = 4 (1 + o_{{\d}}(1))
   \log t_2 - 2 (1 + o_{{\d}}(1)) \log t_1.
\end{equation}
In fact, these yield the inequality
$$
    J_{\rho}(\var_{(\vartheta_1, \vartheta_2)}) \leq (2 \rho_1 - 8 \pi
    + o_\d(1))  \log t_1 + (2 \rho_2 - 8 \pi + o_\d(1)) \log t_2 \to -
    \infty \qquad \quad \hbox{ as } \nu \to 0
$$
uniformly for $(\vartheta_1, \vartheta_2) \in \mathcal{X}_\nu$,
since $\rho_1, \rho_2 > 4 \pi$. Here again we are using that $C
\geq h_i(x) \geq \frac{1}{C}>0$ for any $x \in \Sg$.

We begin by showing \eqref{eq:estlin1}, whose proof clearly also yields \eqref{eq:estlin2}.
It is convenient to write
$$
  \var_1 = \log  \left(1 + \tilde{t}_2^2 d(x_2,y)^2
  \right) - 2 \log \left( 1 + \tilde{t}_1^2 d(x_1,y)^2 \right),
$$
and to divide $\Sg$ into the two subsets
$$
   \mathcal{A}_1 = B_{x_1}(\d) \cup B_{x_2}(\d); \qquad
   \qquad \mathcal{A}_2 = \Sg \setminus \mathcal{A}_1.
$$
For $y \in \mathcal{A}_2$ we have that
$$
  \frac{1}{C_{\d,\Sg} t_1^2} \leq 1 + \tilde{t}_1^2 d(x_1,y)^2 \leq
  \frac{C_{\d,\Sg}}{t_1^2}; \qquad \qquad
   \frac{1}{C_{\d,\Sg} t_2^2} \leq 1 + \tilde{t}_2^2 d(x_2,y)^2 \leq
  \frac{C_{\d,\Sg}}{t_2^2},
$$
which implies
\begin{equation}\label{eq:int1111}
    \frac{1}{|\Sg|} \int_{\mathcal{A}_2} \var_1 dV_g =
   4 (1 + o_{{\d}}(1)) \log t_1 - 2 (1 + o_{{\d}}(1)) \log t_2.
\end{equation}
On the other hand, working in normal geodesic coordinates at $x_i$ one also finds
$$
   \int_{B_{\d}(x_i)} \log \left( 1 + \tilde{t}_i^2 d(x_i,y)^2 \right)
   dV_g = o_\d(1) \log t_i.
$$
Using \eqref{eq:int1111} and the last formula we then obtain \eqref{eq:estlin1}.

Let us now show \eqref{eq:estQ}. We clearly have that
\begin{eqnarray*}
  \n \var_1 & = &
  \n \log \left( 1 + \tilde{t}_2^2
   d(x_2,y)^2 \right) - 2 \n \log \left( 1 + \tilde{t}_1^2 d(x_1,y)^2 \right) \\
  & = & \frac{2 \tilde{t}_2^2 d(x_2,y) \n_y d(x_2,y)}{1 + \tilde{t}_2^2 d(x_2,y)^2}
  - \frac{4 \tilde{t}_1^2 d(x_1,y) \n_y d(x_1,y)}{1 + \tilde{t}_1^2 d(x_1,y)^2},
\end{eqnarray*}
and similarly
\begin{eqnarray*}
  \n \var_2 & = &
  \n \log \left( 1 + \tilde{t}_1^2
   d(x_1,y)^2 \right) - 2 \n \log \left( 1 + \tilde{t}_2^2 d(x_2,y)^2 \right) \\
  & = & \frac{2 \tilde{t}_1^2 d(x_1,y) \n_y d(x_1,y)}{1 + \tilde{t}_1^2 d(x_1,y)^2}
  - \frac{4 \tilde{t}_2^2 d(x_2,y) \n_y d(x_2,y)}{1 + \tilde{t}_2^2 d(x_2,y)^2}.
\end{eqnarray*}
From now on we will assume, without loss of generality, that $t_1 \leq t_2$.
We distinguish between the case $t_2 \geq \d^3$ and $t_2 \leq \d^3$.

In the first case the function $1 + \tilde{t}_2^2 d(x_2,y)^2$ is uniformly
Lipschitz with bounds depending only on $\d$, and therefore we can write that
$$
   \n \var_1 =
   - \frac{4 \tilde{t}_1^2 d(x_1,y) \n_y d(x_1,y)}{1 + \tilde{t}_1^2 d(x_1,y)^2}
   + O_\d(1); \qquad \quad
   \n \var_2 =  \frac{2 \tilde{t}_1^2 d(x_1,y) \n_y d(x_1,y)}{1
   + \tilde{t}_1^2 d(x_1,y)^2} + O_\d(1).
$$
Given a large but fixed constant $C_1 > 0$, we divide the surface $\Sg$ into the three
regions
\begin{equation}\label{eq:3regions}
   \mathcal{B}_1 = B_{x_1}(C_1 t_1); \qquad \quad \mathcal{B}_2 = B_{x_2}(C_1 t_2);
  \qquad \quad \mathcal{B}_3 = \Sg \setminus (\mathcal{B}_1 \cup \mathcal{B}_2).
\end{equation}
In $\mathcal{B}_1$ we have that $|\n \var_i| \leq
{C}{\tilde{t}_1}$, while
\begin{equation}\label{eq:sim1}
    \frac{\tilde{t}_1^2 d(x_1,y) \n_y d(x_1,y)}{1 + \tilde{t}_1^2 d(x_1,y)^2} =
  (1 + o_{C_1}(1)) \frac{ \n_y d(x_1,y)}{d(x_1,y)} \qquad \quad \hbox{ in }
  \Sg \setminus \mathcal{B}_1.
\end{equation}
The last gradient estimates imply that
\begin{eqnarray}\label{eq:estQt2large} \nonumber
  \int_{\Sg} Q(\var_{(\vartheta_1, \vartheta_2)}) dV_g & = &
  \int_{\Sg \setminus \mathcal{B}_1} Q(\var_{(\vartheta_1, \vartheta_2)})
  dV_g + o_{\d}(1) \log \frac{1}{t_1} + O_\d(1)
  \\ & = & 8 \pi \int_{C_1 t_1}^1 \frac{dt}{t} + o_{\d}(1) \log \frac{1}{t_1} + O_\d(1)
 \\ & = & 8 \pi (1 + o_{\d}(1)) \log \frac{1}{t_1} + 8 \pi (1 + o_{\d}(1))
  \log \frac{1}{t_2} + O_\d(1); \qquad \quad t_2 \geq \d^3. \nonumber
\end{eqnarray}
Assume now that $t_2 \leq \d_3$. Then by the definition of $\mathcal{X}_{\nu}$
we have that $d(x_1, x_2) \geq \frac{\d^2}{2}$, and therefore $\mathcal{B}_1
\cap \mathcal{B}_2 = \emptyset$. Similarly to \eqref{eq:sim1} we find
$$
  \left\{
    \begin{array}{ll}
     \frac{\tilde{t}_1^2 d(x_1,y) \n_y d(x_1,y)}{1 + \tilde{t}_1^2
     d(x_1,y)^2} = (1 + o_{C_1}(1)) \frac{ \n_y d(x_1,y)}{d(x_1,y)};  &  \\
     \frac{\tilde{t}_2^2 d(x_2,y) \n_y d(x_2,y)}{1 + \tilde{t}_2^2
     d(x_2,y)^2} = (1 + o_{C_1}(1)) \frac{ \n_y d(x_2,y)}{d(x_2,y)}
      &
    \end{array}
  \right. \qquad \quad \hbox{ in } \mathcal{B}_3.
$$
Moreover we have the estimates
$$
    \left | \n \var_i \right |\leq {C}{\tilde{t}_i} \quad \hbox{ in }
    \mathcal{B}_i, \ i=1,\ 2; \qquad \qquad \left |\n \var_i \right | \leq
    C \quad \hbox{ in } \mathcal{B}_j, \ i \neq j.
$$
Then, there follows:
\begin{eqnarray}\label{eq:estQt2small} \nonumber
  \int_{\Sg} Q(\var_{(\vartheta_1, \vartheta_2)}) dV_g & = &
  \int_{\mathcal{B}_3} Q(\var_{(\vartheta_1, \vartheta_2)}) dV_g +
  o_{\d}(1) \log \frac{1}{t_1} + o_{\d}(1)
  \log \frac{1}{t_2} + O_\d(1) \\
  & = & 8 \pi (1 + o_{\d}(1)) \log \frac{1}{t_1} + 8 \pi (1 + o_{\d}(1))
  \log \frac{1}{t_2} + O_\d(1); \qquad \quad t_2 \leq \d^3.
\end{eqnarray}
With formulas \eqref{eq:estQt2large} and \eqref{eq:estQt2small}, we conclude the proof
of \eqref{eq:estQ} and hence that of the lemma.
\end{pf}

\

\noindent Since the functional $J_\rho$ attains large negative
values on the above test functions $\var_{(\vartheta_1,
\vartheta_2)}$, these are mapped to $X$ by $\Psi$. We next
evaluate the image of $\Psi$ with more precision, beginning with
the following technical lemma.

\begin{lem}\label{l:concscale} Let $\var_1, \var_2$ be as in \eqref{eq:varvar12}:
then, for some $C=C(\d,\Sg)>0$, the following estimates hold
uniformly in $(\vartheta_1, \vartheta_2) \in \mathcal{X}_\nu$:
\begin{equation}\label{eq:noconct1d}
    \sup_{x \in \Sg} \int_{B_x(r t_i)} e^{\var_i} dV_g
   \leq C r^2 \frac{t_i^2}{t_j^2} \qquad \quad \forall r >0,\ i \neq j.
\end{equation}
Moreover, given any $\e > 0$ there exists $C=C(\e, \d, \Sg)$,
depending only on $\e$, $\d$ and $\Sg$ (but not on $\nu$), such
that
\begin{equation}\label{eq:noconct1d3}
  \int_{B_{x_i}(C t_i)} e^{\var_i} dV_g \geq (1 - \e)
  \int_{\Sg} e^{\var_i(y)} dV_g, \ i=1,\ 2.
\end{equation}

uniformly in $(\vartheta_1, \vartheta_2) \in \mathcal{X}_\nu$.
\end{lem}

\begin{pf} We prove the case $i=1$. Observe that
$1 + \tilde{t}_2^2 d(x_2,y)^2 \leq \frac{C}{t_2^2}$ and that $1 +
\tilde{t}_1^2 d(x_1,y)^2 \geq 1$. Therefore we immediately find
$$
  \int_{B_x(t_1 r)} e^{\var_1} dV_g \leq \frac{C}{t_2^2} \int_{B_x(t_1 r)}
  \frac{1}{\left( 1 + \tilde{t}_1^2 d(x_1,y)^2 \right)^2} dV_g(y) \leq C r^2
  \frac{t_1^2}{t_2^2} \qquad \quad \hbox{ for all } x \in \Sg,
$$
which gives the first inequality in \eqref{eq:noconct1d}.

We now show \eqref{eq:noconct1d3}, by evaluating the integral in
the complement of $B_{x_1}(R t_1)$ for some large $R$. Using again
the fact that $1 + \tilde{t}_2^2 d(x_2,y)^2 \leq \frac{C}{t_2^2}$
we clearly have that
\begin{equation}\label{eq:miao}
    \int_{\Sg \setminus B_{x_1}(R t_1)} e^{\var_1(y)} dV_g(y)
  \leq \frac{C}{t_2^2} \int_{\Sg \setminus B_{x_1}(R t_1)}
  \frac{1}{\left( 1 + \tilde{t}_1^2 d(x_1,y)^2 \right)^2} dV_g(y).
\end{equation}
To evaluate the last integral one can use normal geodesic coordinates centered
at $x_1$ and \eqref{eq:intfalt} with a change of variable to find that
$$
   \lim_{t_1 \to 0^+} t_1^{-2} \int_{\Sg \setminus B_{x_1}(R t_1)}
  \frac{1}{\left( 1 + \tilde{t}_1^2 d(x_1,y)^2 \right)^2} dV_g = o_R(1)
  \qquad \quad \hbox{ as } R \to + \infty.
$$
This and \eqref{eq:miao}, jointly with the second inequality in
\eqref{eq:inttot}, conclude the proof of the
\eqref{eq:noconct1d3}, by choosing $R$ sufficiently large,
depending on $\e, \d$ and $\Sg$.
\end{pf}

\

\noindent We next show that, parameterizing the test functions on $\mathcal{X}_{\nu}$
and composing with $R_{\nu} \circ \Psi$, we obtain a map homotopic
to the identity on $\mathcal{X}_{\nu}$. This step will be fundamental for us in
order to run the variational scheme later in this section.

\begin{lem}\label{l:homid} Let $L > 0$ be so large that $\Psi(\{ J_\rho \leq - L \}) \in X$,
and let $\nu$ be so small that $J_\rho(\var_{(\vartheta_1, \vartheta_2)}) < - L$
for $(\vartheta_1, \vartheta_2) \in \mathcal{X}_{\nu}$ (see Lemma \ref{l:dsmallIlow}). Let $R_{\nu}$ be
the retraction given in Lemma \ref{l:retr}. Then the map from $T_\nu : \mathcal{X}_{\nu}
\to \mathcal{X}_{\nu}$ defined as
$$
   T_\nu((\vartheta_1, \vartheta_2)) = R_{\nu} (\Psi(\var_{(\vartheta_1, \vartheta_2)}))
$$
is homotopic to the identity on $\mathcal{X}_{\nu}$.
\end{lem}

\begin{pf} Let us denote $\vartheta_i= (x_i, t_i)$,

$$f_i= \frac{e^{\var_i}}{\int_{\Sg} e^{\var_i} dV_g}, \quad \psi
(f_i)=(\beta_i, \s_i),$$ where $\psi$ is given in Proposition
\ref{covering}. First, we claim that there is a constant
$C=C(\d,\Sg)>0$, depending only on $\Sg$ and $\d$, such that:
\begin{equation}\label{eq:betatest} \frac{1}{C} \leq
\frac{\s_i}{t_i} \leq C,
    \qquad \qquad  d \left( \b_i , x_i
    \right) \leq C t_i.
\end{equation}

By \eqref{eq:noconct1d3}, we have that

$$ \s\left(x_i, f_i \right) \leq C
t_i,$$ where $\s(x,f)$ is the continuous map defined in
\eqref{sigmax}. From that, we get that $\s_i \leq C t_i$. Using
now \eqref{eq:noconct1d}, we get the relation $t_i \leq C \s_i$.

Taking into account that $\s(x_i, f) \leq C t_i$ and \eqref{dett},
we obtain that

$$d \left(x_i, S\left(f_i
\right) \right )\leq C t_i,$$ where $S(f)$ is the set defined in
\eqref{defS}. But since the inequality
$$d\left(\b_i, S\left(f_i
\right) \right) \leq C \s_i$$ is always satisfied, we conclude the
proof of \eqref{eq:betatest}.

\medskip We are now ready to prove the lemma. Let us define a first deformation $H_1$ in the following form:

$$
  \left( \left(
    \begin{array}{c}
      (\b_1, \s_1) \\
      (\b_2, \s_2) \\
    \end{array}
  \right),
s \right) \;\; \stackrel{\small H_1}{\longmapsto} \;\;  \left(
                                \begin{array}{c}
                                  \left( \b_1,\
(1-s) \s_1 + s \kappa_1 \right) \\
                                   \\
\left( \b_2,\ (1-s) \s_2 + s \kappa_2 \right)
                                \end{array}
                              \right),
$$
where $\kappa_i= \min \{ \delta, \frac{\s_i}{\sqrt{\nu}} \}$.

A second deformation $H_2$ is defined in the following way:

$$
  \left( \left(
    \begin{array}{c}
      (\b_1, \kappa_1) \\
      (\b_2, \kappa_2) \\
    \end{array}
  \right),
s \right) \;\; \stackrel{\small H_2}{\longmapsto} \;\;  \left(
                                \begin{array}{c}
                                  \left( (1-s)\b_1 + s x_1, \
\kappa_1 \right) \\
                                   \\
\left( (1-s)\b_2 + s x_2,\ \kappa_2 \right)
                                \end{array}
                              \right),
$$
where $(1-s)\b_i + s x_i$ stands for the geodesic joining
$\beta_i$ and $x_i$ in unit time. A comment is needed here. If
$\kappa_i < \delta$, then we have that $\s_i < \sqrt{\nu} \delta$.
By choosing $\nu$ small enough, this implies that $\beta_i$ and
$x_i$ are close to each other (recall \eqref{eq:betatest}).
Instead, if $\kappa_i= \delta$, the identification in $\ov{\Sg}_\d$ makes the above
deformation trivial.

We also use a third deformation $H_3$:

$$
  \left( \left(
    \begin{array}{c}
      (x_1, \kappa_1) \\
      (x_2, \kappa_2) \\
    \end{array}
  \right),
s \right) \;\; \stackrel{\small H_3}{\longmapsto} \;\;  \left(
                                \begin{array}{c}
                                  \left( x_1,\
(1-s) \kappa_1 + s t_1 \right) \\
                                   \\
\left( x_2,\ (1-s) \kappa_2 + s t_2 \right)
                                \end{array}
                              \right).
$$

We define $H$ as the concatenation of those three homotopies.
Then,

$$ ((\vartheta_1, \vartheta_2), s) \mapsto R_{\nu} \circ H(\Psi(\var_{(\vartheta_1,
\vartheta_2)}),s)$$ gives us the desired homotopy to the identity.

Observe that, since $\nu \ll \d$, $H(\Psi(\var_{(\vartheta_1,
\vartheta_2)}),s)$ always stays in $X$, so that $R_{\nu}$ can be
applied.

\end{pf}

\

\noindent  We now introduce the variational scheme which yields
existence of solutions: this remaining part follows the ideas of
\cite{djlw} (see also \cite{mal}).

Let $\ov{\mathcal{X}}_{\nu}$ denote the (contractible) cone over
$\mathcal{X}_{\nu}$, which can be represented as
$$
   \ov{\mathcal{X}}_{\nu} = \left( \mathcal{X}_{\nu} \times [0,1]
   \right)|_{\sim},
$$
where the equivalence relation $\sim$ identifies $\mathcal{X}_{\nu} \times
\{1\}$ to a single point.
We choose $L > 0$ so large that \eqref{eq:psilow} holds, and then $\nu$ so
small that
$$
  J_{\rho}(\var_{(\vartheta_1, \vartheta_2)}) \leq - 4L \qquad \quad
\hbox{ uniformly for }
   (\vartheta_1, \vartheta_2) \in \mathcal{X}_{\nu},
$$
the last claim being possible by Lemma \ref{l:dsmallIlow}. Fixing this value
of $\nu$, consider the following class of functions
\begin{equation}\label{eq:PiPi}
    \Gamma = \left\{ \eta : \ov{\mathcal{X}}_{\nu} \to H^1(\Sg)
 \; : \; \eta \hbox{ is continuous and } \eta(\cdot \times \{0\})
 = \var_{(\vartheta_1, \vartheta_2)} \hbox{ on } \mathcal{X}_{\nu} \right\}.
\end{equation}
Then we have the following properties.

\begin{lem}\label{l:min-max}
The set $\Gamma$ is non-empty and moreover, letting
$$
  \alpha = \inf_{\eta \in \Gamma}
  \; \sup_{m \in \ov{\mathcal{X}}_{\nu}} J_\rho(\eta(m)), \qquad
  \hbox{ one has } \qquad \alpha > - 2 L.
$$
\end{lem}

\begin{pf}
To prove that $\Gamma \neq \emptyset$, we just notice that the map
\begin{equation}\label{eq:ovPi}
  \tilde{\eta}(\vartheta,s) = s \var_{(\vartheta_1, \vartheta_2)},
  \qquad \qquad (\vartheta,s) \in \ov{\mathcal{X}}_{\nu},
\end{equation}
belongs to $\Gamma$.

Suppose by contradiction that $\alpha \leq  - 2L$: then there
would  exist a map $\eta \in \Gamma$ satisfying the condition
$\sup_{m \in \ov{\mathcal{X}}_{\nu}} J_\rho(\eta(m)) \leq - L$.
Then, since Lemma \ref{l:homid} applies, writing $m = (\vartheta,
s)$, with $\vartheta \in \mathcal{X}_{\nu}$, the map
$$
  s \mapsto R_{\nu} \circ \Psi \circ \eta(\cdot,s)
$$
would be a homotopy in $\mathcal{X}_{\nu}$ between $R_{\nu} \circ
\Psi \circ \var_{(\vartheta_1, \vartheta_2)}$ and a constant map.
But this is impossible since $\mathcal{X}_{\nu}$ is
non-contractible (by the results in Section \ref{s:app} and by the
fact that $\mathcal{X}_{\nu}$ is a retract of $X$) and since
$R_{\nu} \circ \Psi \circ \var_{(\vartheta_1, \vartheta_2)}$ is
homotopic to the identity on $\mathcal{X}_{\nu}$. Therefore we
deduce $\alpha
> - 2 L$, which is the desired conclusion.
\end{pf}

\

From the above Lemma, the functional $J_{\rho}$ satisfies suitable
structural properties for min-max theory. However, we cannot
directly conclude the existence of a critical point, since it is
not known whether the Palais-Smale condition holds or not. The
conclusion needs a different argument, which has been used
intensively (see for instance \cite{djlw, dm}), so we will be
sketchy.

\

\noindent We take $\mu > 0$ such that $\mathcal{J}_i :=
[\rho_i-\mu, \rho_i+\mu]$ is contained in $(4 \pi, 8 \pi)$ for
both $i = 1, 2$. We then consider $\tilde{\rho}_i \in
\mathcal{J}_i$ and the functional $J_{\tilde{\rho}}$ corresponding
to these values of the parameters.

Following the estimates of the previous sections, one easily
checks that the above min-max scheme applies uniformly for
$\tilde{\rho}_i \in \mathcal{J}_i$ for $\nu$ sufficiently small.
More precisely, given any large number $L > 0$, there exists $\nu$
so small that for $\tilde{\rho}_i \in \mathcal{J}_i$
\begin{equation}\label{eq:min-maxrho}
   \sup_{m \in \partial
   \ov{\mathcal{X}}_{\nu}} J_{\tilde{\rho}}(m) < - 4 L; \qquad
   \qquad \alpha_{\tilde{\rho}} := \inf_{\eta \in \Gamma}
  \; \sup_{m \in \ov{\mathcal{X}}_{\nu}} J_{\tilde{\rho}}(\eta(m)) > -
  2L, \ \ (\tilde{\rho}=(\tilde{\rho}_1, \tilde{\rho}_2)).
\end{equation}
where $\Gamma$ is defined in \eqref{eq:PiPi}. Moreover, using for
example the test map \eqref{eq:ovPi}, one shows that for $\mu$
sufficiently small there exists a large constant $\ov{L}$ such
that
\begin{equation}\label{eq:ovlovl}
  \alpha_{\tilde{\rho}} \leq \ov{L} \qquad \qquad
  \hbox{ for  } \tilde{\rho}_i \in \mathcal{J}_i.
\end{equation}

\

\noindent Under these conditions, the following Lemma is
well-known, usually taking the name "monotonicity trick". This
technique was first introduced by Struwe in \cite{struwe}, and
made general in \cite{jeanjean} (see also \cite{djlw, lucia}).

\begin{lem}\label{l:arho}
Let $\nu$ be so small that \eqref{eq:min-maxrho} holds. Then the
functional $J_{t \rho}$ possesses a bounded Palais-Smale sequence
$(u_l)_l$ at level $\tilde{\alpha}_{t \rho}$ for almost every $t
\in \left[ 1 - \frac{\mu}{16 \pi}, 1 + \frac{\mu}{16 \pi}
\right]$.
\end{lem}

\

\begin{pfn} {\sc of Theorem \ref{t:main}.} The existence of a bounded Palais-Smale
sequence for $J_{t \rho}$ implies by standard arguments that this functional possesses
a critical point. Let now $t_j \to 1$, $t_j \in \L$ and let
$(u_{1,j}, u_{2,j})$ denote the corresponding solutions. It is then
sufficient to apply the compactness result in Theorem \ref{th:jlw}, which yields convergence
of $(u_{1,j}, u_{2,j})_j$ by the fact that $\rho_1, \rho_2$ are not multiples of $4 \pi$.
\end{pfn}

\section{Appendix: the set $ X=\ov{\Sg}_{\d} \times \ov{\Sg}_{\d} \setminus
\ov{D}_\d$ is not contractible.}\label{s:app}

Without loss of generality, we consider the case $\delta=1$ (see
\eqref{cono}). Let us denote $\ov{\Sg} = \ov{\Sg}_{1}$.

If $\Sg= \Sf^2$, we have a complete description of $X$. Indeed, in
this case $\ov{\Sg}$ can be identified with $B(0,1) \subset \R^3$.
Therefore, we have:

$$ X = (B(0,1) \times B(0,1)) \setminus E, $$
where $E=\{x \in \R^6:\ x_i = x_{i+3},\ i=1,\ 2,\ 3\}$. By taking
the orthogonal projection onto $E^{\bot}$, we have that $X \simeq
U \setminus \{0\}$ ($\simeq$ stands for homotopical equivalence),
where $U \subset E^{\bot}$ is a convex neighborhood of $0$. And,
clearly, $U\setminus \{0\} \simeq \Sf^2$.

\medskip The case of positive genus is not so easy and we have a less complete
description of $X$. However, we will prove that it is
non-contractible by studying its cohomology groups $H^*(X)$, where
coefficients will be taken in $\R$. Indeed, we will show that:

\begin{pro} \label{pro51} If the genus of $\Sg$ is positive, then $H^4(X)$ is nontrivial. \end{pro}

\begin{pf} In what follows, the elements of $\ov{\Sg}$ will be
written as $(x,t)$, where $x \in \Sg$, $t \in (0,1]$.

Clearly, $X= Y \cup Z$, where $Y$, $Z$ are open sets defined as:

$$Y= \{((x_1,t_1), (x_2,t_2)) \in \ov{\Sg} \times \ov{\Sg}:\ t_1 \neq
t_2 \},$$
$$Z= \{((x_1,t_1), (x_2,t_2)) \in \ov{\Sg} \times \ov{\Sg}:\ t_1<1,\ t_2<1,\
x_1 \neq x_2 \}.$$

Then, the Mayer-Vietoris Theorem gives the exactness of the
sequence:

$$ \cdots \rightarrow H^3(X) \rightarrow H^3(Y) \oplus H^3(Z) \rightarrow H^3(Y \cap Z) \rightarrow H^4(X) \rightarrow \cdots  $$

Since our coefficients are real, the above cohomology groups are
indeed real vector spaces. The exactness of the sequence then gives:
\begin{equation} \label{mv}  dim(H^3(Y \cap Z)) \leq dim(H^4(X)) + dim(H^3(Y) \oplus H^3(Z)).\end{equation}

Let us describe the sets involved above. First of all, observe
that $Y=Y_1 \cup Y_2$ has two connected components:
$$Y_i= \{((x_1,t_1), (x_2,t_2)) \in \ov{\Sg} \times \ov{\Sg}:\ t_i
> t_j, \ j \neq i \}.$$

To study $Y_1$, we define the following deformation retraction:

$$ r_1 : Y_1 \to Y_1, \qquad r_1((x_1,t_1), (x_2,t_2))= ((x_1,1),
(x_2,1/2)).$$

Clearly, $r_1(Y_1)=0 \times \left(\Sg \times \{1/2\}\right)$,
which is homeomorphic to $\Sg$.
Analogously, $Y_2 \simeq \Sg$, and so $Y \simeq \Sg \dotcup \Sg$
(here $\dotcup$ stands for the disjoint union).

For what concerns $Z$, we can define a deformation retraction:

$$r : Z \to Z, \qquad r((x_1,t_1), (x_2,t_2))= ((x_1,1/2),
(x_2,1/2)).$$

Observe that $r(Z)= \left( \Sg \times \{ 1/2\} \times \Sg \times
\{ 1/2\} \right)\setminus \ov{D}$ which is homeomorphic to $\Sg
\times \Sg \setminus D$, where $D$ is the diagonal of $\Sg \times
\Sg$. Let us set
$$
A=\Sg \times \Sg \setminus D,
$$
since it will appear many times in what follows.

Moreover, $Y \cap Z = (Y_1 \cap Z) \cup (Y_2 \cap Z)$, and so this
has two connected components. Also here we have a deformation
retraction:

$$ r'_1: Y_1 \cap Z \to Y_1 \cap Z ,\qquad r'_1((x_1,t_1), (x_2,t_2))= ((x_1,1/2),
(x_2,1/3)).$$

It is clear that $r'_1(Y_1 \cap Z)$ is homeomorphic to $A=\Sg
\times \Sg \setminus D$. Analogously we can argue for $Y_2 \cap
Z$; therefore, $Y \cap Z \simeq A \dotcup A$.

Hence, from \eqref{mv} we obtain:
\begin{equation} \label{H30} dim(H^4(X)) \geq dim (H^3(A)).\end{equation}

\medskip Let us now compute the cohomology  of $A=\Sg \times
\Sg \setminus D$. Given $\e>0$, let us define:
$$B = \{(x,y) \in \Sg \times \Sg:\ d(x,y) < \e  \},$$ which is an open neighborhood of $D$. Clearly,
we can use the local contractibility of $\Sg$ to retract $B$ onto
$D$. Moreover, $A \cup B = \Sg \times \Sg$. The Mayer-Vietoris
Theorem yields the exact sequence:

\begin{equation} \label{mv2}\cdots \rightarrow H^2(A \cap B) \rightarrow H^3(\Sg \times \Sg)
 \rightarrow H^3(A) \oplus H^3(B) \rightarrow H^3(A\cap B) \rightarrow \cdots \end{equation}

Therefore, in order to study $H^3(A)$ we need some information
about $H^*(A\cap B)$.

By using the exponential map, we can define a homeomorphism:

$$ h : A \cap B = \{(x,y) \in \Sg \times \Sg:\ 0< d(x,y) < \e  \} \to \{(x,v) \in T\Sg:\ 0<\|v\|<\e\},$$
$$ h(x,y)= (x,v) \in T\Sg \mbox{ such that } \exp_x(v)=y,$$
where $T \Sg$ is the tangent bundle of $\Sg$. Therefore, $A \cap
B$ is homotopically equivalent to the unit tangent bundle $UT
\Sg$.

The cohomology of $UT\Sg$ must be well known, but we have not been
able to find a precise reference. We state and prove the following
lemma:

\begin{lem} Let us denote by $g=g(\Sg)$ the genus of $\Sg$. Then:

\begin{enumerate}

\item if $g= 1$, $H^0(UT \Sg)\cong H^3(UT \Sg)\cong \R$ and
$H^1(UT \Sg)\cong H^2(UT \Sg)\cong \R^3$.

\item if $g \neq 1$, $H^0(UT \Sg) \cong H^3(UT \Sg)\cong \R$ and
$H^1(UT \Sg)\cong H^2(UT \Sg)\cong \R^{2g}$.

\end{enumerate}

\end{lem}

\begin{pf} We only need to compute $H^1(UT \Sg)$ and $H^2(UT
\Sg)$.

If $g=1$, that is, $\Sg \simeq \T^2$, then $T\Sg$ is trivial and
hence $UT \Sg \simeq \T^2 \times \Sf^1 \simeq \T^3$. The K{\"u}nneth
formula gives us the result.

If $g \neq 1$, we use the Gysin exact sequence (see Proposition
14.33 of \cite{bott-tu}):

$$0 \rightarrow H^1(\Sg) \rightarrow H^1(UT\Sg) \rightarrow
H^0(\Sg) \stackrel{\wedge e}{\rightarrow} H^2(\Sg) \rightarrow
H^2(UT\Sg) \rightarrow H^1(\Sg) \rightarrow H^3(\Sg)=0.$$

In the above sequence, $\wedge e$ is the wedge product with the
Euler class $e$. Since we are working with real coefficients and
the Euler characteristic of $\Sg$ is different from zero, then
$\wedge e$ is an isomorphism. Therefore, we conclude:
$$  H^1(UT\Sg) \cong H^1(\Sg) \cong \R^{2g}, \qquad  H^2(UT\Sg) \cong H^1(\Sg) \cong \R^{2g}.$$

\end{pf}

\begin{rem} We have chosen real coefficients since they simplify
our arguments and are enough for our purposes. As a counterpart,
the above computations do not take into account the torsion part.
For instance, it is known that $UT \Sf^2 = \R \mathbb{P}^3$ (see
\cite{montesinos}).
\end{rem}

We now come back to the proof of Proposition \ref{pro51}. With our
information, \eqref{mv2} becomes:

$$ \cdots \rightarrow H^2(A \cap B) \rightarrow \R^{4g}
 \rightarrow H^3(A) \rightarrow \R \rightarrow \cdots
 $$
In the above sequence we  computed $H^3(\Sg \times \Sg)$
using the K{\"u}nneth formula. Then, $4g \leq dim(H^2(UT \Sg)) +
dim(H^3(A))$. Therefore, $dim(H^3(A)) \geq 2g$, if $g>1$, or
$dim(H^3(A)) \geq 1$, if $g=1$. In any case we conclude by
\eqref{H30}.

\end{pf}


\begin{thebibliography}{99}


%
%
%
%
%
%



\bibitem{bott-tu}{R. Bott and L. W. Tu, }{Differential forms in algebraic
topology, }{Graduate Texts in Mathematics, 82. Springer-Verlag,
New York-Berlin, 1982.}


%
%


%
%
%
%
%
%
%
%

%
%



\bibitem{cl}{W. X. Chen and C. Li, }{ Prescribing Gaussian curvatures on
surfaces with conical singularities, }{J. Geom. Anal. 1-4, (1991),
359-372.}

\bibitem{cl2}{W. X. Chen and C. Li, }{Classification of solutions of some nonlinear elliptic equations, }
{Duke Math. J. 63 (1991), no. 3, 615-622.}

\bibitem{clin} C.C Chen and C.S. Lin, { Topological degree for a mean
field equation on Riemann surfaces}, Comm. Pure Appl. Math. 56-12
(2003), 1667-1727.

\bibitem{sw}{M. Chipot, I. Shafrir and G. Wolansky, }{On the solutions of Liouville systems, }
{J. Differential Equations 140 (1997), no. 1, 59-105.}

\bibitem{demarchis}{F. De Marchis, }{Generic multiplicity for a scalar field equation
on compact surfaces, }{J. Funct. Anal. 259 (2010) 2165-2192.}


\bibitem{djlw}{W. Ding, J. Jost, J. Li and G. Wang, }{Existence results
for mean field equations, }{Ann. Inst. Henri Poincar\'e, Anal. Non
Lin{\`e}aire 16-5 (1999), 653-666.}


\bibitem{dja}{Z. Djadli, }{ Existence result for the mean field problem
on Riemann surfaces of all genus}, Comm. Contemp. Math. 10 (2008), no. 2, 205-220.






\bibitem{dm}{Z. Djadli and A. Malchiodi, }{Existence of
conformal metrics with constant $Q$-curvature, }{Annals of Math.,
168 (2008),  no. 3, 813-858.}


\bibitem{dunne}{G. Dunne, }{Self-dual Chern-Simons Theories, }{Lecture
Notes in Physics, vol. 36, Berlin: Springer-Verlag, 1995.}

\bibitem{guest}{M. A. Guest, }{Harmonic maps, loops groups, and integrable systems. }{London Mathematical
Society Student Texts, 38. Cambridge University Press, Cambridge,
1997.}

%
\bibitem{hat}{A. Hatcher, }{Algebraic Topology, }{Cambridge
University Press 2002.}

\bibitem{jeanjean}{L. Jeanjean, }{On the existence of bounded Palais-Smale sequences and applications
 to a Landesman-Lazer type problem set on $R^N$, }{Proc. Roy.
Soc. Edinburgh A 129 (1999) 787-809.}



%

\bibitem{jlw}{J. Jost, C. S. Lin and G. Wang, }{Analytic aspects of the Toda system II.
Bubbling behavior and existence of solutions, }{Comm. Pure Appl.
Math. 59 (2006), 526-558.}

\bibitem{jw2}{J. Jost and G. Wang, }{Classification of solutions of a Toda
system in $\R^2$, }{Int. Math. Res. Not., 2002 (2002), 277-290.}

\bibitem{jw}{J. Jost and G. Wang, }{Analytic aspects of the Toda system I. A Moser-Trudinger
inequality, }{Comm. Pure Appl. Math. 54 (2001), 1289-1319.}



%
%

%


\bibitem{ll}{J. Li and Y. Li, }{Solutions for Toda systems on Riemann surfaces, }{Ann. Sc. Norm. Super. Pisa Cl. Sci.
(5) 4 (2005), no. 4, 703-728.}

\bibitem{lin-zhang}{C.S. Lin and L. Zhang, }{A Topological Degree Counting for Some Liouville
Systems of Mean Field Type, }{Comm. Pure Appl. Math. 64 (2011),
556-590.}

\bibitem{lucia}{M. Lucia, }{A deformation lemma with an application to
a mean field equation, }{ Topol. Methods Nonlinear Anal.  30
(2007), no. 1, 113-138.}

%
%
%


\bibitem{mreview}{A. Malchiodi, }{Topological methods for
an elliptic equation with exponential nonlinearities, }{Discrete
Contin. Dyn. Syst.  21  (2008),  no. 1, 277-294.}



\bibitem{mal}{A. Malchiodi, }{Morse theory and a scalar
field equation on compact surfaces, }{Adv. Diff. Eq., 13 (2008),
1109-1129.}

\bibitem{cheikh}{A. Malchiodi and C. B. Ndiaye, }{Some existence results for the Toda system on closed surfaces, }
{Atti Accad. Naz. Lincei Cl. Sci. Fis. Mat. Natur. Rend. Lincei
(9) Mat. Appl. 18 (2007), no. 4, 391-412.}

\bibitem{mr}{A. Malchiodi and D. Ruiz, }{ New improved Moser-Trudinger
inequalities and singular Liouville equations on compact surfaces,
}{to appear in GAFA, preprint arXiv:1007.3861.}


\bibitem{montesinos}{J. M. Montesinos, }{Classical Tesselations
and Three-Manifolds, }{Springer-Verlag Berlin-Heidelberg, 1987.}

%
%


%
%




\bibitem{sw2}{I. Shafrir and G. Wolansky, }{Moser-Trudinger and logarithmic HLS inequalities for systems, }{J. Eur.
Math. Soc. 7-4 (2005), 413-448.}


\bibitem{struwe}{M. Struwe, }{On the evolution of harmonic mappings of Riemannian surfaces, }{Comment.
Math. Helv. 60 (1985) 558-581.}





%
%


\bibitem{st}{M. Struwe and G. Tarantello, }{ On multivortex solutions in
Chern-Simons gauge theory}, {Boll. Unione Mat. Ital., Sez. B,
Artic. Ric. Mat. (8)-1, (1998), 109-121.}

\bibitem{tar}{G. Tarantello, }{ Self-Dual Gauge Field Vortices:
An Analytical Approach, }{PNLDE 72, Birkh\"auser Boston, Inc.,
Boston, MA, 2007.}

\bibitem{tar3}{G. Tarantello, }{Analytical, geometrical and topological
aspects of a class of mean field equations on surfaces, }{Discrete
Contin. Dyn. Syst.  28  (2010),  no. 3, 931-973.}

\bibitem{wang}{G. Wang, }{Moser-Trudinger inequalities and Liouville systems, }
{C. R. Acad. Sci. Paris. S{\'e}r. I Math. 328 (1999), no. 10,
895-900.}

\bibitem{yys}{Y. Yang, }{Solitons in Field Theory and Nonlinear
Analysis, }{Springer-Verlag, 2001.}



\end{thebibliography}
\end{document}